\documentclass[10pt]{article}
\usepackage{geometry}                
\usepackage{graphicx}
\usepackage{subfig}
\usepackage{amssymb}
\usepackage{amsmath}
\usepackage{hyperref}
\usepackage{epstopdf}
\usepackage{caption}
\newcommand{\vecb}{\boldsymbol{\mathbf{b}}}

\newcommand{\vecw}{\boldsymbol{\mathbf{w}}}
\newcommand{\vecx}{\boldsymbol{\mathbf{x}}}

\newcommand{\matA}{\boldsymbol{\mathbf{A}}}
\newcommand{\matB}{\boldsymbol{\mathbf{B}}}

\newcommand{\matF}{\boldsymbol{\mathbf{F}}}

\newcommand{\matzero}{\boldsymbol{\mathbf{0}}}

\newcommand{\tv}{\boldsymbol{\theta}}
\newcommand{\mv}{\boldsymbol{\mu}}

\newcommand{\Sv}{\boldsymbol{\Sigma}}

\title{Information-adaptive clinical trials with\\selective recruitment and binary outcomes}
\author{James E. Barrett\footnote{Contact: james.barrett@kcl.ac.uk}\\
\small{Division of Health \& Social Care Research, King's College London}\\
\small{UCL Cancer Institute, University College London}}
\date{30 May, 2017}

\begin{document}
\maketitle

\begin{abstract}

Selective recruitment designs preferentially recruit individuals that are estimated to be statistically informative onto a clinical trial. Individuals that are expected to contribute less information have a lower probability of recruitment. Furthermore, in an information-adaptive design recruits are allocated to treatment arms in a manner that maximises information gain. The informativeness of an individual depends on their covariate (or biomarker) values and how information is defined is a critical element of information-adaptive designs. In this paper we define and evaluate four different methods for quantifying statistical information. Using both experimental data and numerical simulations we show that selective recruitment designs can offer a substantial increase in statistical power compared to randomised designs. In trials without selective recruitment we find that allocating individuals to treatment arms according to information-adaptive protocols also leads to an increase in statistical power. Consequently, selective recruitment designs can potentially achieve successful trials using fewer recruits thereby offering economic and ethical advantages.
\end{abstract}

%
%
\section{Introduction}
\label{sec:intro}
%
%

Selective recruitment designs have recently been proposed in which patients that are expected to provide more statistical information are more likely to be recruited onto a clinical trial, otherwise they are rejected \cite{Barrett16}. Informativeness depends on a patient's covariate values and the aim is to avoid wasting resources on covariate values that are uninformative and instead focus on regions of the covariate space where we expect to learn more. We use the term covariates to mean both clinical covariates of interest or biomarkers. Since the cohort is enriched with the ``most informative'' patients fewer recruits are required overall thereby leading to potential economic savings and ethical benefits.

The principal aim in a clinical trial is to establish a statistical relationship between treatments, covariates and clinical outcomes. Selective recruitment designs are motivated by the observation that not all patients provide the same amount of statistical information towards this goal. For instance, a patient with covariate values that are identical to patients already in the cohort may not be as helpful as a patient with previously unseen covariate values. Alternatively, if there is a part of the covariate space where the association with clinical outcomes is ambiguous then it may be beneficial to recruit more patients from within that region in order to resolve the ambiguity. Finally, in some circumstances extreme covariate values have a more pronounced relationship with outcomes and enrichment of those covariates may help to elucidate how outcomes depend on covariates. In all of these cases, it may be beneficial to selectively target individuals that are expected to be informative. In an information-adaptive design regardless of whether individuals are selected for informativeness or not the probability of allocation to a treatment arm is proportional to the amount of information they are expected to contribute on that arm. Information-adaptive designs are adaptive in the sense that what is deemed to be informative will depend on the observations accrued so far in the trial.

To illustrate the advantages of a selective recruitment design, consider a trial that aims to establish whether a biomarker is predictive of treatment. Suppose measuring the biomarker costs \$10 and recruitment costs \$100 per patient (which would include the cost of drugs, outcome measurement, administrative costs etc.). Compare a successful randomised trial which requires 100 patients (at a total cost of $100\times\$100 + 100\times\$10 = \$11,000$) to a selective recruitment design in which 50 patients out of 200 candidates are recruited (at a cost of $50\times\$100 + 200\times\$10 = \$7,000$). We will see shortly that these cohort sizes are quite typical. Aside from the obvious economic savings an ethical argument can be made in favour of exposing fewer patients to treatments with uncertain efficacies (which may include placebos). Additionally, patients that are rejected from the trial are freed for recruitment onto other trials instead which may benefit research areas with limited eligible patients.

A key element of selective recruitment designs is to define a \emph{utility} function that will quantify how \emph{useful} or \emph{informative} an individual is (as a function of their covariates). In \cite{Barrett16} selective recruitment designs were applied to time-to-event outcomes and informativeness was based on the posterior entropy of a proportional hazards model. In this paper we extend the concept of selective recruitment to trials with binary outcomes and now consider four different ways of measuring informativeness.

The first is a heuristic approach which we refer to as \emph{uncertainty sampling}. Patients for which either binary outcome is equally likely are targeted since by concentrating on patients about whom we are most uncertain we aim to achieve a classifier that will generalise to previously unseen patients better and achieve better predictive performance in future. The second method, the \emph{posterior entropy} method, utilises Shannon's entropy which provides a measure of uncertainty regarding a random variable by selecting samples that are expected to minimise the posterior entropy \cite{MACK92}. Thirdly, selecting patients that minimise the expected \emph{generalisation error} (the proportion of incorrect classifications) places more emphasis on predictive performance than parameter uncertainty. Finally, the \emph{variance reduction} strategy uses Fisher's information matrix in order to estimate the variance of predictions and selects patients that are expected to minimise this. We develop some general theory for how these four methods can be used in clinical trials with binary outcomes. We then implement it in the specific case where a logistic regression model is assumed to connect the binary outcomes to covariates.

Several of the concepts above have previously been applied to the field of machine learning under the term \emph{query learning} (where we wish to query the underlying system in an optimal way) or \emph{active learning} (where we actively seek out informative data samples). In the statistics literature this is known as \emph{optimal experimental design}. In \cite{schein2007active} several active learning techniques were applied to logistic regression and a thorough evaluation is provided using simulated datasets and experimental data. The authors used a \emph{pool-based} sampling approach where a pool of unlabelled samples exists and the task is to select which samples should be selected for labelling. A thorough overview of active learning in general can be found in \cite{settles2010active}. Active learning approaches have been used to determine \emph{individualised treatment rules} by selectively recruiting informative patients in \cite{minsker2015, deng2011}. Covariates are included in \cite{minsker2015} and the design is similar in spirt to the uncertainty sampling approach since patients for whom the optimal treatment was ambiguous were recruited. These approaches are applied to two treatments and continuous outcomes.

In the theory of optimal experimental design Fisher's information matrix is useful because it provides a lower bound to the variance of our parameter estimator (this is known as the Cram\'{e}r-Rao inequality). \emph{$D$-optimal} designs maximise the determinant of the information matrix. We shall see that in our application to logistic regression this will be equivalent to the posterior entropy approach. \emph{$A$-optimal} designs minimise the trace of the inverse information matrix and in our case this is equivalent to the variance reduction method.

An information-adaptive design can be conducted with or without selective recruitment. Even if all candidates are accepted onto the trial they can still be allocated in a way that maximises statistical information gain. There are numerous adaptive clinical trial designs that allocate patients in a manner that aims to maximise either the total benefit to recruited patients or the statistical information gained. \emph{Response-adaptive} designs will increase the probability with which patients are allocated to treatment arms that appear most efficacious. \emph{Covariate-adaptive} designs try to ensure a balanced distribution of covariates across treatment arms making it easier to compare the treatment effects \cite{pocock1975sequential, taves1974minimization}. \emph{Covariate adjusted response-adaptive} designs will select the best treatment arm after taking covariate values into account \cite{rosenberger2001covariate}. This obviously benefits individuals recruited onto the trial but it may not lead to a statistically optimal study. Other proposals are based on optimising explicit statistical quantities. A $D$-optimal adaptive design based on the theory of optimal experimental design has previously been developed by \cite{atkinson1982optimum}. The emphasis in this paper is on statistical considerations only and we will not examine response-adaptive protocols. A good overview of adaptive designs can be found in the textbook \cite{yin2013clinical}. 

In Section \ref{eq:modeldef} we define the four measures of informativeness and describe the selective recruitment protocol. Results from an experimental dataset are presented in Section \ref{sec:casestudies} and results from numerical simulations are given in Section \ref{sec:numerics}. Finally, some of the practical aspects of selective recruitment designs are discussed in Section \ref{sec:disc}.

%
%
\section{Information-adaptive clinical trial designs}
\label{eq:modeldef}
%
%

\subsection{Model definition}

Observed data are $D_n = \{(\vecx_1,y_1),\ldots,(\vecx_{n},y_{n})\}$, and consist of pairs of covariate vectors $\vecx\in\mathbb{R}^d$ and binary labels $y\in\{-1,+1\}$. The number of current recruits is $n$ (and increases as the trial progresses) and $N$ is the total number of recruits. We assume there are $K$ treatment arms. The relationship between outcomes and covariates on arm $k$ is completely specified by $p(y|\vecx,\tv_k)$ where $\tv_k$ is a vector of model parameters. Some elements of $\tv_k$ may be common across some (or all) of the arms, but for simplicity we will assume that this is not the case here. The posterior over model parameters for arm $k$ is given by Bayes' rule
\begin{equation}
p(\tv_k|D_n) = \frac{p(D_n|\tv_k)p(\tv_k)}{p(D_n)},
\end{equation}
where the marginal likelihood is $p(D_n) = \int \text{d} \tv_k \,p(D_n|\tv_k)p(\tv_k)$. For i.i.d. samples the data likelihood term is $p(D_n|\tv_k) = \prod_{i\in R_k} p(y_i|\vecx_i,\tv_k)$ where $R_k$ is the subset of patients allocated to arm $k$. Letting $\left<\cdots\right>_p$ denote the expectation with respect to $p$, predictions for an individual on arm $k$ with covariates $\vecx$ are given by $p_k(y|\vecx,D_n) = \left<p(y|\vecx,\tv_k)\right>_{p(\tv_k|D_n)}$.
%
\subsection{Utility functions}
%

We now define precisely our four utility functions. These will be later used to calculate recruitment and arm allocation probabilities for a candidate patient with covariates $\vecx^*$. The utility function is denoted by $E_k(\vecx^*|D_n)$ and measures how informative $\vecx^*$ is considered to be if they were allocated to arm $k$, with larger values corresponding to more informative covariates.

\subsubsection{Uncertainty sampling.}
Here samples with uncertain predictions are regarded as informative. The quantity $1-p(\hat{y}|\vecx^*,\tv_k)$ where $\hat{y} = \text{argmax}_y p(y|\vecx^*,\tv_k)$ provides a measure of our uncertainty in predicting the class corresponding to $\vecx^*$. The utility function is given by
\begin{equation}
E_k(\vecx^*|D_n) = \left<1-p(\hat{y}|\vecx^*,\tv_k)\right>_{p(\tv_k|D_n)}.
\label{eq:01loss}
\end{equation}

\subsubsection{Posterior entropy.}

The entropy of the posterior distribution for arm $k$ is $S_{k}(D_n) = -\left<\log p(\tv_k|D_n)\right>_{p(\tv_k|D_n)}$ by definition. The expected entropy in the hypothetical scenario where the candidate $(\vecx^*,y^*)$ is added to arm $k$ is $S_k(\vecx^*|D_n) = \left<S_k(D_n\cup(\vecx^*,y^*))\right>_{p_k(y^*|\vecx^*,D_n)}$ where we take the expectation with respect to the unknown outcome $y^*$ (using the predictive distribution) and $D_n\cup(\vecx^*,y^*)$ is the union of the current dataset and the candidate datum. The utility function is given by the decrease in the expected posterior entropy
\begin{equation}
E_k(\vecx^*|D_n) = S_k(D_n) - S_k(\vecx^*|D_n).
\label{eq:postent}
\end{equation}

\subsubsection{Generalisation error.}

The expected generalisation error for arm $k$ is obtained by taking the expectation of (\ref{eq:01loss}) with respect to all possible inputs and defining $\epsilon_k(D_n) = \big\langle\left<1-p(\hat{y}|\vecx,\tv_k)\right>_{p(\tv_k|D_n)}\big\rangle_{p(\vecx)}$. The expected variance if candidate $(\vecx^*,y^*)$ were to be added to arm $k$ is $\epsilon_k(\vecx^*|D_n) = \left<\epsilon_k(D_n\cup(\vecx^*,y^*))\right>_{p_k(y^*|\vecx^*,D_n)}$ and the utility function is defined as
\begin{equation}
E_k(\vecx^*|D_n) = \epsilon_k(D_n) - \epsilon_k(\vecx^*|D_n).
\end{equation}

\subsubsection{Variance reduction.}

The covariance matrix of $\tv$ about the posterior mode $\hat{\tv}$ (for conciseness we drop the arm index $k$) is approximated by Fisher's information matrix $\matF$ with $F_{\rho\nu} = -\partial^2/\partial \theta_{\rho} \partial \theta_{\nu}\log p(\tv|D_n) |_{\tv=\hat{\tv}}$. The variance of an arbitrary function $g(\tv)$ can be approximated by $\nabla g\cdot\matF^{-1}\nabla g$ where $[\nabla g]_{\nu}=\partial g(\tv)/\partial \theta_{\nu}$ \cite{MACK92}. We will consider the variance of the predictive density and define $\tilde{\nu}^2(\vecx|D_n) = \nabla p(y|\vecx,\tv)\cdot\matF^{-1}\nabla p(y|\vecx,\tv)$. We then take the expectation over all possible inputs and define $\tilde{\sigma}^2(D_n) = \left<\tilde{\nu}^2(\vecx|D_n)\right>_{p(\vecx)}$. The expected variance if candidate $(\vecx^*,y^*)$ is added to arm $k$ is $\tilde{\sigma}^2(\vecx^*|D_n) = \left<\tilde{\sigma}^2(D_n\cup(\vecx^*,y^*)\right>_{p_k(y^*|\vecx^*,D_n)}$. A utility function is given by the decrease in expected predictive variance
\begin{equation}
E_k(\vecx^*|D_n) = \tilde{\sigma}^2(D_n) - \tilde{\sigma}^2(\vecx^*|D_n).
\end{equation}

\subsubsection{$D$-Optimal and $A$-optimal designs.}

In our specific implementation we use a logistic regression model with a Gaussian approximation of the posterior distribution. We show in Section A.2.2 of the Supplementary Information that in this case the $D$-optimality criterion is equivalent to the posterior entropy approach since both methods use the determinant of $\matF$ as a measure of informativeness. We also show in Section A.2.4 of the Supplementary Information that under the assumption of Gaussian distributed covariates the variance reduction method is equivalent to an $A$-optimality criterion. By selecting data points to optimise either the determinant or trace of $\matF$ we generally end up with a narrower and more sharply defined posterior distribution. The Wald test statistic, which we use to test for significant model parameters, essentially captures this feature and consequently it can reject the null hypothesis using fewer observations.

\subsection{Allocation and recruitment rules}

The utility function is used to calculate both a treatment allocation and recruitment  probability via the quantity
\begin{equation}
\rho_k(\vecx^*|D_n) = \frac{E_k(\vecx^*|D_n) - E_{min}^k(D_n)}{E^k_{max}(D_n) - E^k_{min}(D_n)}.
\label{eq:rho}
\end{equation}
where $E^k_{max}(D_n) = \max_{\vecx}E_k(\vecx|D_n)$ and $E^k_{min}(D_n)=\min_{\vecx}E_k(\vecx|D_n)$. The quantity $\rho(\vecx^*|D_n)$ provides a measure of where $\vecx^*$ falls between the maximally and minimally informative covariates and takes a value of 1 (or 0) when $\vecx^*$ is maximally (or minimally) informative. A treatment arm is drawn from
\begin{equation}
p(k|\vecx^*,D_n) = \frac{\rho_k(\vecx^*|D_n)}{\sum_{j=1}^K  \rho_j(\vecx^*|D_n)}\quad\text{for $k=1,\ldots,K$.}
\label{eq:biased}
\end{equation}
Having selected a treatment arm $k^*$ a recruitment probability is simply given by $\rho_{k^*}(\vecx^*|D_n)\in[0,1]$. Therefore, candidates with higher utility are more likely to be recruited.

It is possible to conduct a trial without selective-recruitment but with adaptive allocation. A treatment arm is assigned at random according to (\ref{eq:biased}) and the recruitment probability is simply set to one. In this case no patients are rejected, but statistical power may be gained due to the informative allocation of patients to treatment arms.

\begin{figure}[tb]
\centering
\includegraphics[scale=0.62]{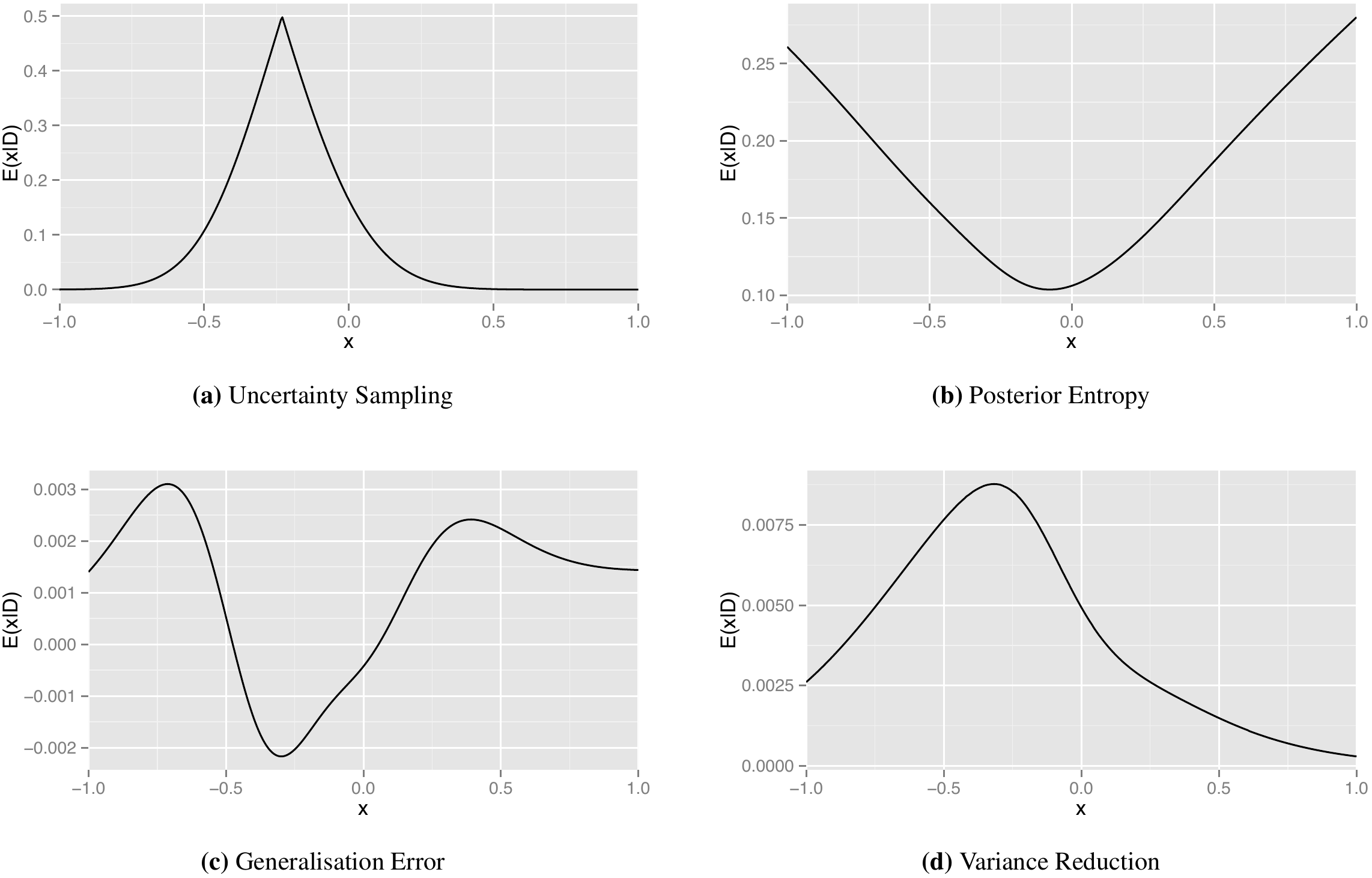}
\caption{Case Study 1: Plots of the four utility functions versus $x^*$. The data observed so far consist of five burn-in patients from the Wisconsin Diagnostic Breast Cancer dataset.}
\label{fig:case1}
\end{figure}

Note that in practice the optimal utility values are located within a restricted search space. This is because for some utility functions $E_k(\vecx^*|D_n)\to\pm\infty$ for $\vecx^*\to\pm\infty$ (see Figure \ref{fig:case1} (b) for example). This is problematic since it means that a candidate with $\vecx^*$ is compared to atypical or unrealistic optimal utility values. We found that searching over the hypercube defined by the first and ninth deciles of the marginal population distributions works well in practice. In other words, the hypothetical patients which correspond to the optimal utility values should be patients that are likely to be observed in practice. One consequence of this is that if $\vecx^*$ lies outside of the search space then it is possible that (\ref{eq:rho}) can take values outside $[0,1]$. This is dealt with by simply mapping any values of (\ref{eq:rho}) which are less than zero to zero (and values exceeding one are mapped to one).

Secondly, note that alternative allocation and recruitment strategies can be implemented. Randomised allocation to treatment arms can be achieved with $p(k|\vecx^*,D_n) = 1/K$ for $k=1,\ldots,K$. Deterministic allocation would allocate to $\text{argmax}_{k} \rho_{k}(\vecx^*|D_n)$ with probability one. Alternative recruitment protocols can be achieved by passing $\rho_k(\vecx^*|D_n)$ as an argument to a transformation function $f_0(s)$. For instance $f_0(s)=\theta(s-p_0)$ will deterministically recruit patients that exceed a threshold $p_0$. The \emph{step function} (or heaviside function) is defined by $\theta(s)=1$ for $s>0$ and 0 otherwise. Alternatively, $f_0(s) = (1+\text{tanh}(s/\beta_0 + p_0))/2$ allows one to control the degree of stringency in the recruitment process by tuning the values of $\beta_0$ and $p_0$. When $\beta_0\to0$ this corresponds to deterministic recruitment.

Finally, in practice, it is desirable to impose a \emph{burn-in} period at the very beginning of the trial where a certain number of patients are recruited before selection begins. This helps ensure that the posterior and predictive densities include some observed data before selective recruitment begins.

\subsection{Application to logistic regression}
\label{sec:lr-def}

We apply the above methods to a logistic regression model where for arm $k$ we write $p(y=+1|\vecx,\vecw_k,w_{k0}) = 1/(1+e^{-\vecw_k\cdot\vecx - w_{k0}})$. The vector of model parameters $\tv_k$ is therefore composed of the regression weights $\vecw_k$ and the intercept term $w_{k0}$. All of the four methods above depend on integrals that are analytically intractable. We use various approximations in order to achieve a practical implementation. Full details are given in the Supplementary Information. The prior over model parameters $p(w_0,\vecw)$ is assumed to factorise over each parameter. Gaussian priors are chosen with zero mean and variance $\alpha_0^2$. Unless otherwise stated we used $\alpha_0^2=5$ throughout this paper.

%
%
\section{Case Study: The Wisconsin Diagnostic Breast Cancer dataset}
\label{sec:casestudies}
%
%

\begin{table}[b!]
\centering
\begin{tabular}{{|c|c|c|c|c}}
\hline
Method & Statistical power & Validation success & Rejected\\
\hline
Randomised Trial & 46.4\% & 68.9\% & 0  \\
Uncertainty Sampling & 28.0\% & 68.9\% & 44.9  \\
Posterior Entropy & 81.0\% & 69.4\% & 30.0 \\
Generalisation Error & 65.4\%& 69.1\%& 33.5\\
Variance Reduction & 60.0\% & 69.1\% & 26.0\\
\hline
\end{tabular}
\caption{Case Study: The second column gives the statistical power for $w$ (the logistic regression weight parameter defined in Section \ref{sec:lr-def}). The validation success is the percentage of correct predictions made on a validation dataset. The final column is the average number of rejections. Results obtained by averaging over 500 simulations.}
\label{tab:case1average}
\end{table}

We begin by comparing the performance of the four selective recruitment protocols on a simulated randomised trial based on the Wisconsin Diagnostic Breast Cancer dataset which was downloaded from the UCI Machine Learning Repository. The dataset contains 30 real-valued covariates for a total of 569 patients. Each patient is classified into a benign or malignant group. We choose one of those covariates, tumour `smoothness', for a univariate analysis. For the purposes of simulation we use the tumour smoothness as a predictive biomarker and the benign and malignant groups as an arbitrary binary outcome (in practice, this could be treatment responders or non-responders for instance). The simulated trial contains a single arm in which a binary outcome is measured and our goal is to establish any statistical association between the biomarker and the outcome.

The tumour smoothness variable was transformed linearly such that all values lay in the range $[-1,+1]$. The search space for locating optima of utility functions was defined by the first and ninth decile of the covariate distribution ($-0.8$ and 0.8 respectively). A burn-in period where the first five patients were automatically recruited was imposed. A total of 25 patients were recruited onto the trial. Since the order in which individuals arrive was not recorded in the original dataset we permuted the arrival order to simulate different trials. Note that the same order of patients is used when comparing the different methods.

In order to gain some insight into how the four different utility functions work we have plotted them as a function of $x^*$ in Figure \ref{fig:case1}. The functions are plotted immediately after the $n=5$ burn-in individuals were recruited in one particular simulation. The decision boundary is defined as the values of $x$ where $p(y|x,w,w_0)=0.5$. Fitting a logistic regression model using these five observations we find the decision boundary at $x=-0.23$. The uncertainty sampling utility function, in (a), takes a maximum value of $0.5$ at this location as expected, and drops towards zero elsewhere. The characteristic shape of the posterior entropy function is shown in (b). This function tends to have higher values for large absolute values of $x$. This is because terms with large $x$ have a bigger impact on the posterior and consequently reduce the entropy more.

The generalisation error utility function in (c) has a minimum close to the decision boundary at $x=-0.23$ (where $p(y|x,w,w_0)=0.5$). Two local maxima occur at either side of this. Thus, samples that are close --- but not too close --- to the point of highest predictive uncertainty at $x=-0.23$ are favoured. The variance reduction utility function has a shape that is similar to the uncertainty sampling utility but is not as narrowly focused on $x=-0.23$. We therefore expect the variance reduction method to recruit a more diverse cohort than uncertainty sampling. 

The next patient in this particular trial had $x^*=-0.32$. In the case of uncertainty sampling, for example, $E(\vecx^*|D_n) = 0.34$. According to (\ref{eq:rho}) with probabilistic recruitment there is a 68\% chance of recruitment. If we were using the posterior entropy method this probability would be 14\% and so forth.

In Table \ref{tab:case1average} are the results from each method averaged over 500 simulations (each with a different permutation of the order in which patients arrive). The randomised clinical trial (RCT) design achieved a statistical power of 46.4\%. This was calculated as the percentage of simulated trials in which the null hypothesis is correctly rejected. A Wald test was used which compares the inferred $w$ to the null hypothesis of $w=0$. Using the entropy utility function statistical power rose to 81\% although this required an average of 30 rejections. The uncertainty sampling method performed worse than random with only 28\% power. The remaining two utility measures also increased the statistical power above an RCT.

In order to see if selective recruitment led to better predictive accuracy we reserved 25 patients as a validation dataset. The binary clinical outcome was predicted from the biomarker for each individual in the validation cohort (using a logistic regression model fitted to the data acquired during the simulated trial). For each simulated trial a different set of 25 individuals were randomly selected to form the validation cohort. We report that all five methods achieved a similar level of predictive accuracy with only a marginal improvement when using an adaptive design.

%
%
\section{Numerical Simulations}
\label{sec:numerics}
%
%

\subsection{Simulation Study 1: comparison of different covariate distributions}

\begin{figure}[tb!]
\centering
\begin{tabular}{c c}
\subfloat{\includegraphics[scale=0.6]{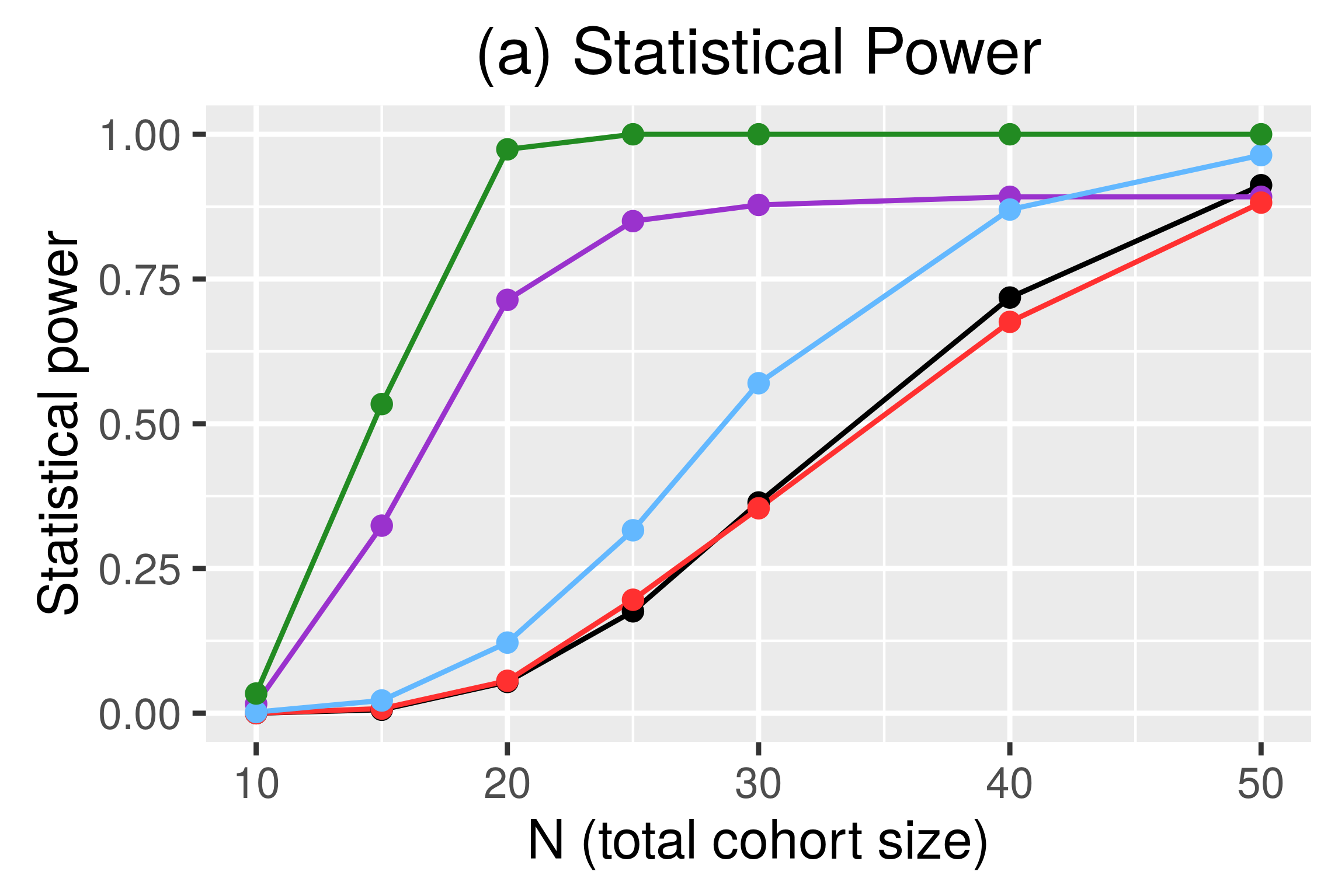}} & \subfloat{\includegraphics[scale = 0.6]{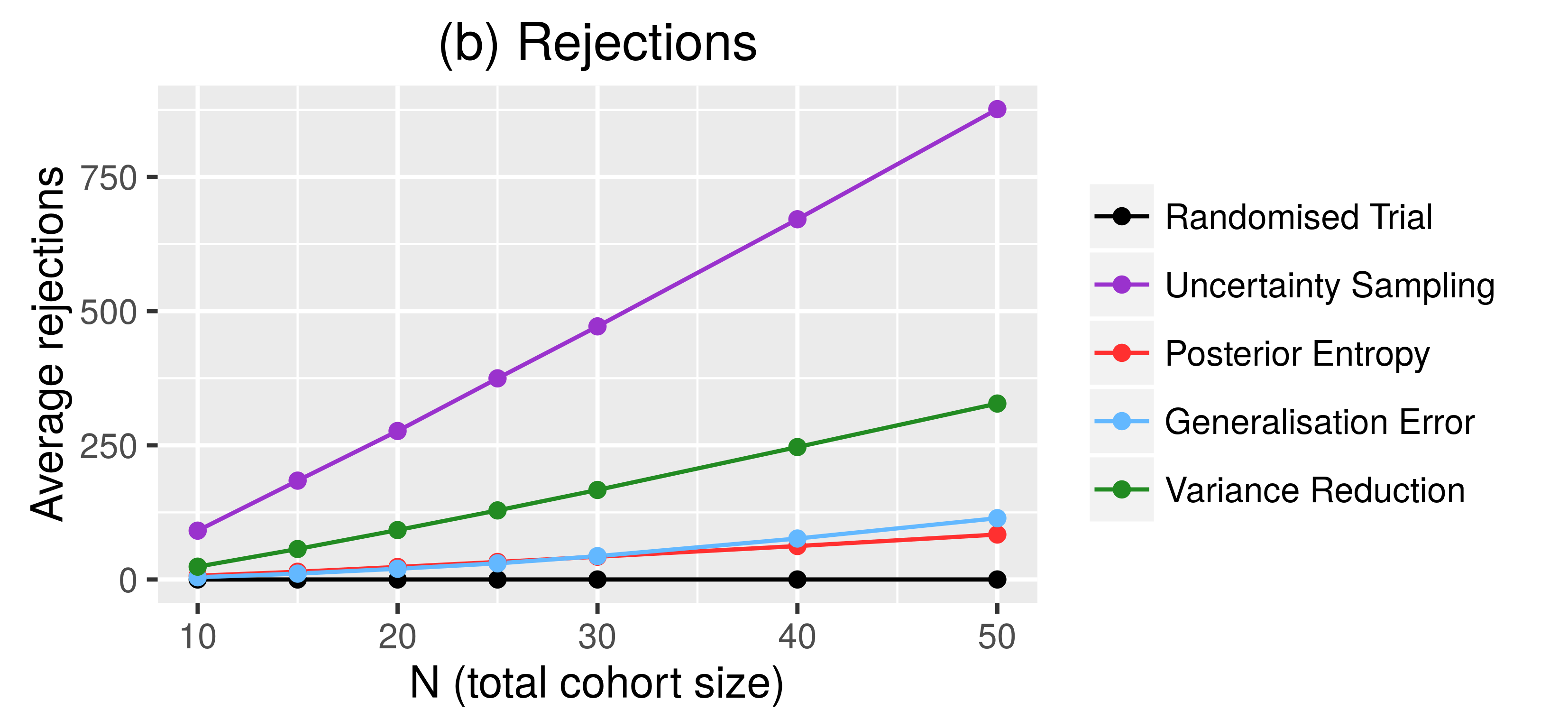}}\\
\subfloat{\includegraphics[scale=0.6]{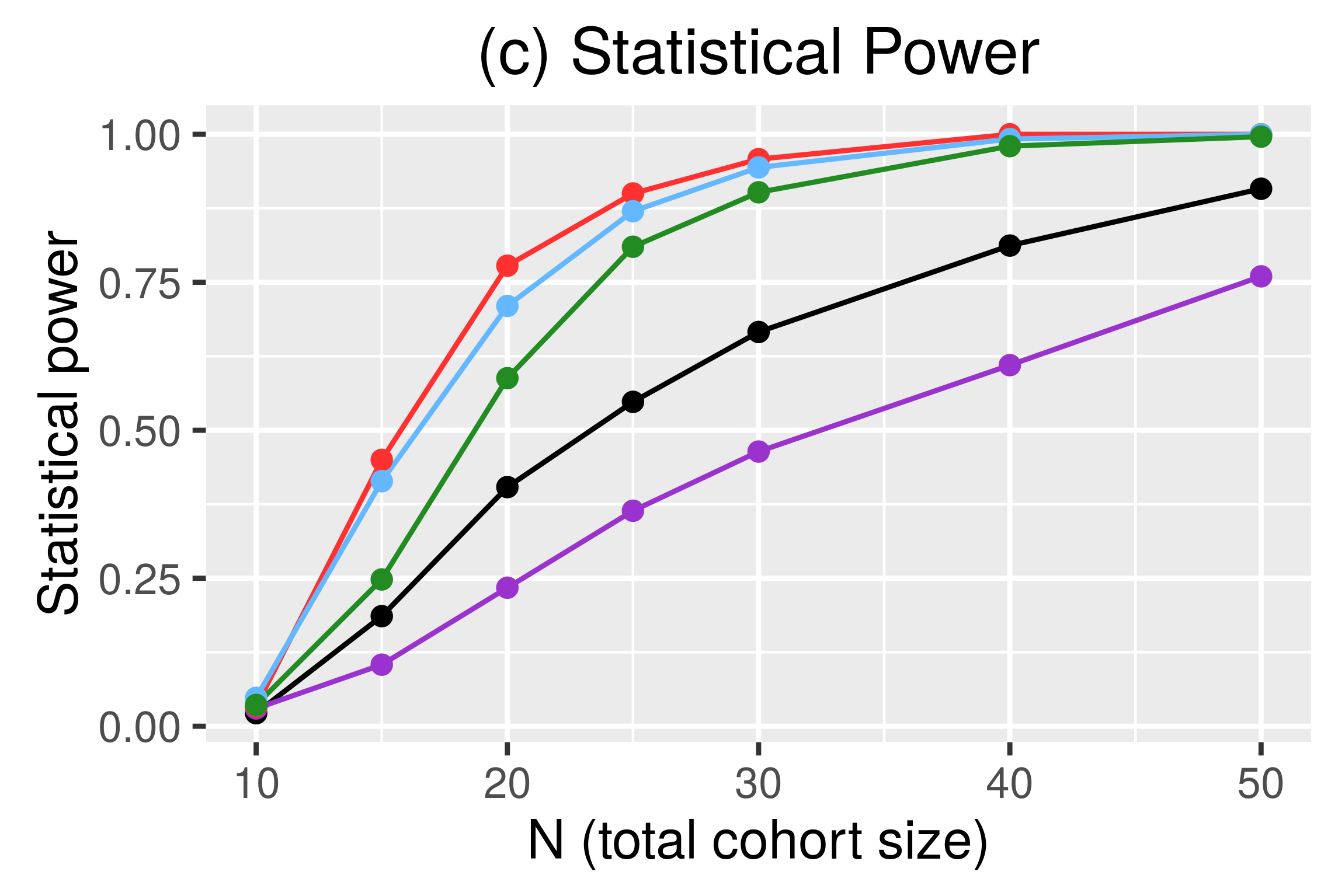}} & \subfloat{\includegraphics[scale = 0.6]{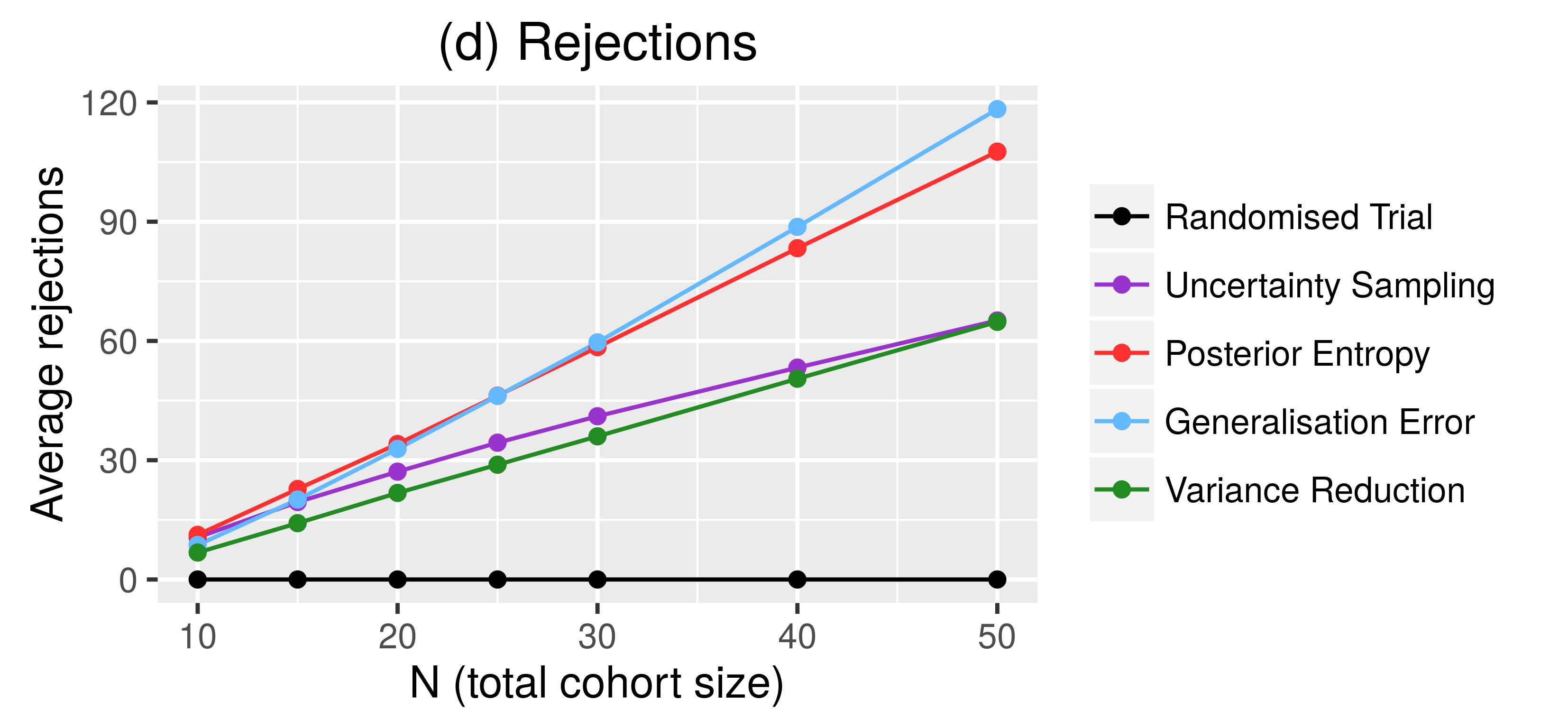}}\\
\end{tabular}
\caption{Simulation Study 1: in (a) and (c) is the statistical power for the weight parameter $w$ and in (b) and (d) is the average number of rejections versus the total cohort size $N$. The top two panels correspond to the uniformly distributed covariates that are almost linearly separable. The bottom two panels correspond to the non-sparable example.}
\label{fig:sim1:1}
\end{figure}

To understand the conditions under which each utility function performs well we simulate trials using two different types of patient distributions that are designed to highlight the strengths and weaknesses of the different methods. In the first case, the \emph{linearly separable} case, one dimensional covariates were generated from a uniform distribution between $\pm1$. Binary labels were generated according to a logistic regression model with $w=32$ and $w_0=-8$. These relatively large parameter values generate classes that are almost perfectly linearly separable with a clearly defined decision boundary at $x=0.25$ (Supplementary Figure 3a). Secondly, in the \emph{non-separable case}, patients with $y=+1$ and $y=-1$ were drawn from two Gaussian distributions with standard deviation 0.5 centred on -0.25 and 0.25 respectively (Supplementary Figure 4a). This represents a more heterogeneous cohort with substantial overlap of both classes. These two extreme cases are intended to highlight the strengths and weaknesses of the four methods and to offer some insight into how they work rather than represent realistic datasets. We perform a more realistic simultion in Section \ref{sec:sim2}. In both cases 500 trials were simulated with a total of 50 patients recruited in each trial. The statistical power and average number of rejections as a function of $N$ were recorded and are plotted in Figure \ref{fig:sim1:1}. Results from validation cohorts are shown in Supplementary Figure 2.

The uncertainty sampling method performs dramatically better than a randomised trial in the linearly separable case (Figure \ref{fig:sim1:1} (a)), although this increase in power comes at the cost of rejecting a comparatively large number of patients (Figure \ref{fig:sim1:1} (b)). Examination of the empirical covariate distribution (Supplementary Figure 3c) shows that only patients that are very close to the decision boundary are recruited. This is beneficial in the linearly separable case since this is the most informative part of the covariate space. It also explains the large number of rejections since only patients very close to $x=0.25$ are recruited. In the non-separable case the uncertainty sampling method performs worse than random because it focuses on the most heterogenous part of the population around $x=0$ and therefore fails to take advantage of information that could be acquired towards $x=\pm1$. This extremely narrow focus on the decision boundary is detrimental.

Conversely, the posterior entropy method works well in the non-separable case because it favours samples towards the extremes of the covariate space --- a part of the covariate space where a clearer picture of the association between covariates and outcomes is revealed. In the linearly separable case this strategy is only marginally more effective than a randomised trial. Although the method does sample some patients close to the decision boundary (Supplementary Figure 3d) it predominantly concentrates resources into the less informative regions at the extremes of the covariate space.

The generalisation error samples to either side of the decision boundary in the linearly separable case although it also appears to sample roughly uniformly from the entire population. In the non-separable case it operates remarkably similarly to the posterior entropy approach and archives reasonably good performance.

The variance reduction method is more versatile and achieves good performance in both cases. Like uncertainty sampling it also concentrates on the decision boundary in the linearly separable case, albeit not as tightly. It therefore rejects fewer patients while obtaining greater statistical power. In the non-separable case it samples patients on either side of the decision boundary (Supplementary Figure 4e), thus avoiding the the most heterogeneous regions.

Note that examination of the cohort distribution the Case Study (Supplementary Figure 1) reveals a distribution that is very similar to the non-separable case considered here. This helps to explain why the uncertainty sampling method performed poorly. Also note that with the linearly separable data we relaxed the priors over model parameters by setting $\alpha_0=20$ due to the large parameter values. This goes some way to explaining why the statistical power is so low in the randomised trial. Larger values of the parameters are a priori more likely and therefore more evidence is required in order to establish statistical significance.

\subsection{Simulation Study 2: a three arm trial with selective recruitment}
\label{sec:sim2}

\begin{figure}[t]
\centering
\begin{tabular}{c c}
\subfloat{\includegraphics[scale=0.6]{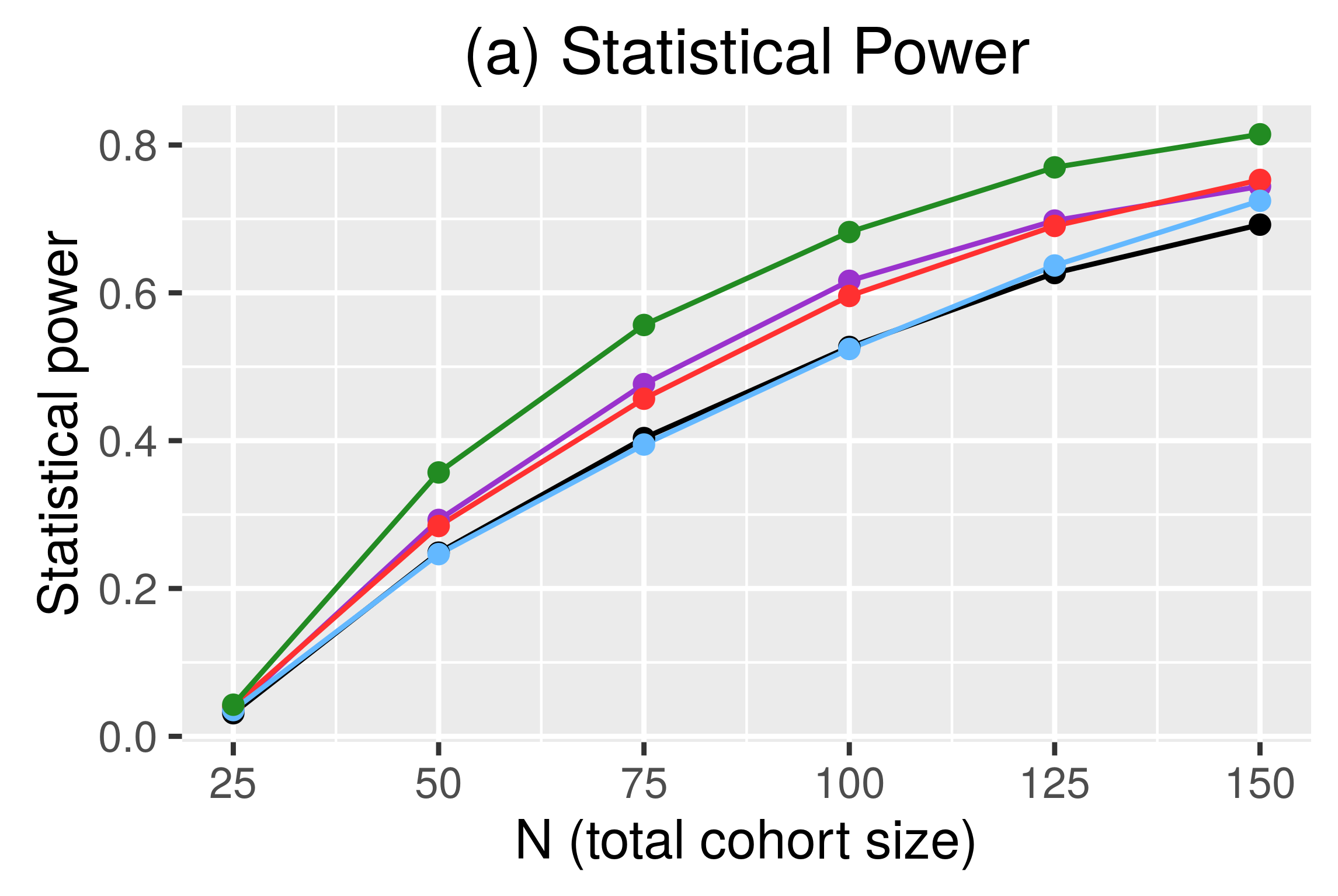}} & \subfloat{\includegraphics[scale = 0.6]{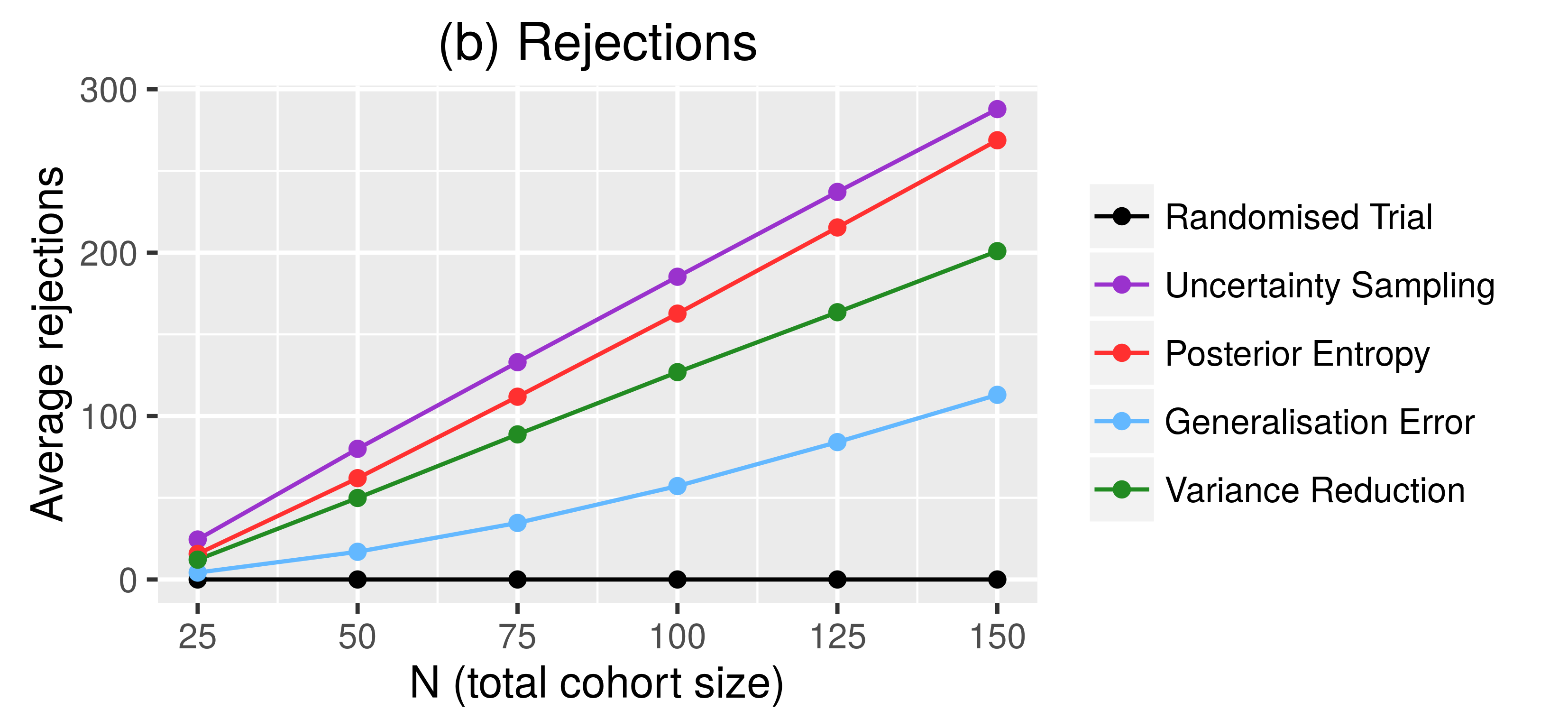}}\\
\end{tabular}
\caption{Simulation Study 2: (a) Statistical power and, (b), average number of rejections versus the number of recruited patients $N$. Results from 500 simulations with three treatment arms, selective recruitment, and information-adaptive treatment allocation.}

\label{fig:sim2:1}
\end{figure}

Here we simulated a more complex trial with three arms and two covariates. Parameter values for each arm of the logistic regression model were set to $\vecw_1 = (-3, 6)$, $\vecw_2 = (4, -8)$, and $\vecw_3=(5, 2)$ and $(w_{01},w_{02},w_{03})=(1.5,-1.5,0)$. Patients were both selectively recruited and allocated to treatment arms adaptively. Covariates were drawn from a uniform distribution over the unit square with vertices at $(\pm1,\pm1)$. A total of 500 trials were simulated. A burn-in period with fifteen patients was imposed. 

In Figure \ref{fig:sim2:1} the statistical power and number of rejections are plotted. The adaptive designs offer consistently higher power than randomised trials with the variance reduction method performing best. An increase in power of up to 16\% is achieved. Examination of Figure \ref{fig:sim2:1} (a) shows that a variance reduction adaptive design with approximately 100 recruits and 125 rejections can achieve the same statistical power as a randomised trial with 150 recruits. Plots of the mean square error between inferred and true parameter values and the validation success rate are shown in Supplementary Figure 5.

We then repeated the simulation with all values of $\vecw$ and $w_0$ equal to zero in order to calculate the Type I error rate. We found that, with the exception of uncertainty sampling, all methods had an error rate of about 5\% (Supplementary Figure 6a). The variance reduction method had a slightly inflated error rate that increases with $N$ for smaller cohort sizes. All methods converged to the true parameter estimates (Supplementary Figure 6b). Uncertainty sampling converged relatively slowly in comparison to the other methods and had a substantially elevated Type I error rate.

\subsection{Simulation Study 3: a three arm trial without selective recruitment}

\begin{figure}[tb]
\centering
\includegraphics[scale=0.65]{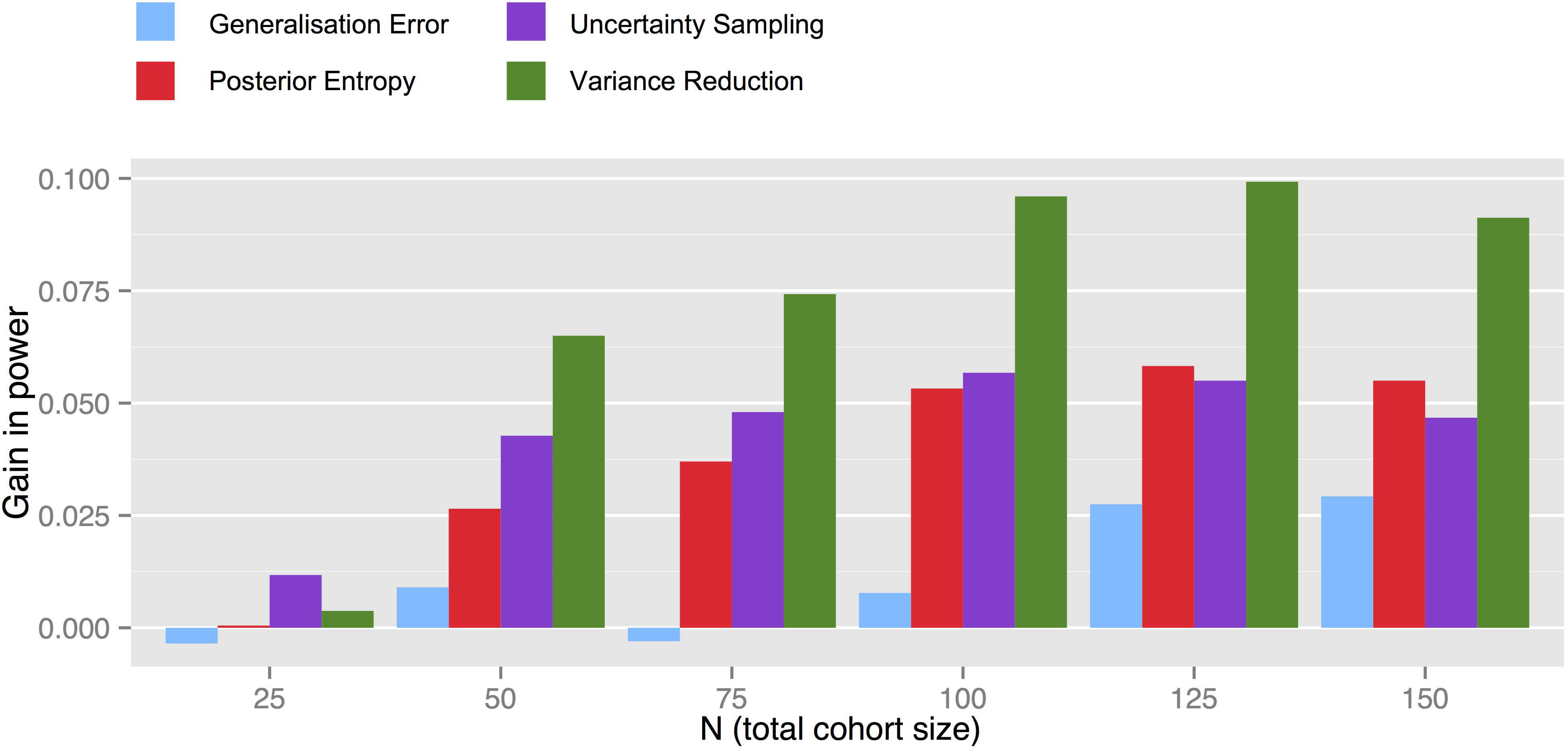}
\caption{Comparison of the absolute gain in statistical power when using information-adaptive allocation and random allocation in Simulation Study 3 (where selective recruitment is not used).}
\label{fig:sim2:2}
\end{figure}

Finally, using the same setup as above we simulate trials with three treatment arms and adaptive allocation but without any selective recruitment. In Figure \ref{fig:sim2:2} the absolute gain in statistical power (compared to the randomised trial) is plotted for each method. Gains of up to 10\% were observed with the variance reduction method without rejecting any individuals. We conclude that adaptive allocation to treatments is beneficial regardless of whether participants are recruited selectively or not although this gain in power is not quite as high as the selective recruitment version above. In supplementary figure 7 we directly compare designs with and without selective recruitment and with and without adaptive allocation.

In the original selective recruitment study of time-to-event outcomes based on the posterior entropy \cite{Barrett16} it was reported that information-adaptive allocation lead to no discernible gain in statistical power. It is now clear that this was due to a subtle difference in how the allocation probabilities were calculated. In the original study $\rho_k(\vecx^*|D_n) = E_k(\vecx^*|D_n)/E^k_{max}(D_n)$ instead of (\ref{eq:rho}). This effectively allocated patients to treatment arms that had fewer recruits. To see why this is the case first note that a treatment arm with fewer patients will generally see a larger decrease in entropy when another patient is added. Therefore both $E^k_{min}(D_n)$ and $E^K_{max}(D_n)$ are larger in arms with fewer patients. Suppose for instance $E_k(\vecx^*|D_n)=4.5$, $E^k_{min}(D_n)=4$ and $E_{max}(D_n)=5$. According to \cite{Barrett16} $\rho_k(\vecx^*|D_n)=4.5/5=0.9$ whereas equation (\ref{eq:rho}) gives $\rho_k(\vecx^*)=0.5$. In the original study arms with fewer patients therefore tended to have inflated allocation probabilities that acted merely to balance treatment arms. The formulation of allocation probabilities in this paper is more intuitive and directs patients to optimal treatment arms.

Finally, we checked if the distribution of patients across arms was significantly imbalanced due to information-adaptive allocation. Since there is a total of 150 recruits per trial randomised allocation would result in an average of 50 patients per arm. We performed a chi squared test for each simulated trial to see if the distribution of patients differed significantly from random. We adjusted for multiple hypotheses testing because we are testing 500 simulations in total. We found that the posterior entropy method deviated from random with statistical significance in 6 simulations (at $p=0.05$). The smallest arm had a median of 43 patients and the biggest arm had a median of 57. In 268 simulations the variance reduction method significantly deviated from random. The smallest arm had a median of 38 patients, the largest had a median of 64. This suggests that despite offering the highest relative gain in power the variance reduction method has a strong tendency to skew the treatment arm sizes which may be undesirable in practice. The generalisation error did not differ from random in any trials. The uncertainty sampling method reached significance in 416 trials with the smallest arm having a median of 32 patients and the largest having a median of 72.

%
%
\section{Discussion}
\label{sec:disc}
%
%

Selective recruitment designs potentially offer greater statistical power than conventional randomised designs with cohorts of the same size. Conversely, a desired level of statistical power can be achieved with fewer recruited patients. Even in trials without selective recruitment it can still be beneficial to allocate individuals to treatment arms according to information-adaptive protocols. Results suggest that statistical power can be substantially boosted by preferentially accruing informative individuals. It is essential that the selection is done prospectively on the basis of covariates only.

The methodology depends on prior knowledge of $p(y|x,D)$. This could be in the form of prior clinical knowledge that the relationship between covariates and outcomes is monotonic and that therefore a logistic regression model (or a proportional hazards model in the case of time-to-event outcomes) is appropriate. Otherwise a sufficient burn-in phase can be used to check model assumptions before proceeding with a selective recruitment protocol.

Of the four different methods of quantifying statistical informativeness considered here we found that the variance reduction method offered consistently good performance under all types of conditions considered here and it tends to reject a comparatively moderate number of patients albeit at the cost of a slightly inflated type I error rate for small cohort sizes. The uncertainty sampling method tends to converge slowly to the correct parameter estimates, has an inflated type I error, rejects a comparatively large number of patients and leads to imbalanced treatment arms, and is therefore not recommended for future use. It has previously been reported that the uncertainty sampling method performs poorly in active machine learning tasks \cite{schein2007active} and our findings are consistent with this observation. The posterior entropy and generalisation error approaches work well for non-separable populations but poorly with data that are almost linearly separable.

There are several practical issues worth noting and selective recruitment designs are not suitable for all types of trials. Firstly, we have considered covariates that are continuously distributed across the population. Populations with more diversity in terms of their covariates, and hence their informativeness, possess greater scope for forming informative cohorts. Such designs may be applicable to categorical covariates but this has not been explored in this paper. Trials conducted on populations with a limited range of covariate values might not benefit from selective recruitment. Secondly, the recruitment rate will be lower if only a subset of all eligible candidates are recruited. Selective recruitment designs will be suitable for scenarios with sufficiently high accrual rates or where longer trial times are acceptable.

A third issue is that the cohort distribution will tend to deviate from the population distribution. Numerical results suggest that this depends strongly on the method used and the nature of the population distribution. This could pose problems for generalising the study findings to the larger population. A simple remedy, however, is to impose a minimum recruitment probability to ensure that there is a sufficient degree of random sampling from the population. On top of this randomly sampled background selective recruitment can enrich the cohort for informative patients to a desired extent. In conclusion, the population distribution and accrual rate will need to be considered a priori in order to determine if a selective recruitment design is appropriate.

Information-adaptation with selective recruitment protocols allow the acquisition of medical evidence using fewer patients than traditional randomised designs. Clinical trials are not intended to be therapeutic and exposing fewer patients to the effects of unproven treatments could offer ethical advantages in some cases. Moreover, patients rejected from one trial are free to participate in alternative trials thereby allowing patients to contribute more effectively to medical research. Reduced costs due to smaller trials means that resources can be diverted towards further research instead. Future investigations will focus on incorporating response-adaptive protocols, in which the treatment response of individual participants is taken into account, in order to establish a more general adaptive framework for answering clinical questions in an efficient and flexible manner.

\section{Acknowledgements}

This research was supported by the CRUK \& EPSRC Comprehensive Cancer Imaging Centre at King's College London and University College London jointly funded by Cancer Research UK and the Engineering and Physical Sciences Research Council (EPSRC). The research was part funded/supported by the National Institute for Health Research (NIHR) Biomedical Research Centre based at Guy's and St Thomas' NHS Foundation Trust and King's College London. The views expressed are those of the author(s) and not necessarily those of the NHS, the NIHR or the Department of Health. The author thanks Louise Browne, Aylin Cakiroglu, Stephan Beck, Tony Ng, and Ton Coolen for helpful discussions. The author acknowledges the use of the UCL Legion High Performance Computing Facility (Legion@UCL), and associated support services in the completion of this work.

\bibliographystyle{unsrt}
\bibliography{refs}

%
%
\newpage
\section*{Supplementary Information}
%
%
\setcounter{figure}{0}

This section contains the practical details of implementing the different methods of quantifying statistical information when a logistic regression model is assumed. Several additional results and plots are also presented here.

%
%
\appendix
\section{Implementation}
\label{eq:implementation}
%
%

We assume a logistic regression model with a variational approximation so that the posterior takes a convenient form. We give details of how to apply the four utility functions described in the main text to this specific scenario.

\subsection{Logistic Regression with a Variational Approximation}
\label{sec:variational}
In a logistic regression model we specify
\begin{equation}
p(y=+1|\vecx,\vecw,w_0) = \frac{1}{1+e^{-\vecw\cdot\vecx - w_0}}
\end{equation}
where $\vecw\in\mathbb{R}^d$ is a vector of regression coefficients and $w_0\in\mathbb{R}$ is the intercept term. The probability of belonging to the $-1$ class is $1 - p(y=+1|\vecx,\vecw)$. For compactness we write the $d+1$ dimensional vector $(w_0,\vecw)$ as $\vecw$. Since the first component of $\vecw$ is now $w_0$ we redefine $\vecx$ as a $d+1$ vector with the first component equal to 1. If there are multiple treatment arms then each arm will have a distinct set of parameters $\vecw_k$ but for simplicity we will drop the index $k$ here. In order to have a more convenient form for the posterior we will use a variational approximation. A detailed description of variational logistic regression is provided in \cite{bishop2006pattern} which we will use but not reproduce here. The posterior is approximated by $p(\vecw|D_n) \approx \mathcal{N}(\mv_n,\Sv_n)$ with
\begin{equation}
\Sv_n^{-1} = \Sv_0^{-1} + 2\sum_{i=1}^{n}\lambda(\xi_i)\vecx_i\cdot\vecx_i^{\text{T}}\quad\text{ and }\quad
\mv_n = \Sv_n\cdot\left(\sum_{i=1}^{n}\frac{y_i}{2}\vecx_i\right)
\end{equation}
where $\Sv_0$ is the diagonal matrix with elements $\alpha_0^2$ from the prior $p(\vecw)$ and
\begin{equation}
\lambda(\xi) = \frac{1}{2\xi}\left(\frac{1}{1+e^{-\xi}} - \frac{1}{2}\right)
\end{equation}
with real valued variational parameters $\xi_1,\ldots,\xi_N$ that can be determined using the EM algorithm as described in \cite[Section 10.6]{bishop2006pattern}. Finally, following the example of \cite[Section 4.5]{bishop2006pattern} we can write the predictive distribution as
\begin{equation}
p(y^*=+1|\vecx^*,D_n) = \frac{1}{1+\exp(-\mv_n\cdot\vecx^*(1+\pi\vecx^*\cdot\Sv_n\vecx^*/8)^{-1/2})}.
\label{eq:approxpred}
\end{equation}

\subsection{Utility Functions}

In this section we discuss the specific implementation of the four utility functions considered in the main text.

\subsubsection{Uncertainty Sampling}

We use the approximated predictive distribution (\ref{eq:approxpred}) to write $E_k(\vecx^*|D_n) = 1 - p_k(\hat{y}|\vecx^*,D_n)$. Moreover, the minimum and maximum values of $E_k$ are, by definition, 0 and 0.5, respectively which means that numerical optimisation is not required. Uncertainty sampling is therefore the fastest of the four methods.

\subsubsection{Posterior Entropy}

Under the variational approximation the posterior takes the form of a multivariate Gaussian. The posterior entropy is given by $S(D_n) = d(1+\log(2\pi))/2 + |\Sv_n|/2$. The Fisher Information matrix is $\matF=\Sv_n^{-1}$ under the variational approximation. Minimisation of the entropy is equivalent to maximising $|\matF|$ as in a $D$-optimal design. To compute the expected entropy if candidate $\vecx^*$ were to be recruited the predictive distribution (\ref{eq:approxpred}) is used and two distinct variational approximations are computed for $S(D_n\cup(\vecx^*,\pm1))$.

\subsubsection{Generalisation Error}

We assume that $p(\vecx)$ is a uniform distribution over the hypercube with edges at $\pm1$ (this is arbitrary so we assume covariates are scaled appropriately). The expected generalisation error is obtained by numerically integrating $1 - p(\hat{y}|\vecx,D_n)$ over $p(\vecx)$. Again, the predictive distribution (\ref{eq:approxpred}) is used.

\subsubsection{Variance Reduction}
\label{sec:imp:vard}

Here we consider the case where $d=2$. We can obtain an analytical solution to the integral $\tilde{\sigma}^2(D_n) = \left<\tilde{\nu}^2(\vecx|D_n)\right>_{p(\vecx)}$ if we approximate the logistic sigmoidal function by a probit function: $1/(1+e^{-z})\approx \Phi(\lambda z)$ with $\lambda^2 = \pi/8$ \cite[Section 4.5.2]{bishop2006pattern}. This gives us
\begin{equation}
\tilde{\sigma}^2(D_n) = \sum_{\mu,\nu=0}^d \frac{\lambda^2}{2\pi}F^{-1}_{\mu\nu}\int \text{d}\vecx\,p(\vecx)x_{\mu}x_{\nu}e^{-\lambda^2(\vecw\cdot\vecx)^2}.
\end{equation}
We can define $\matA\in\mathbb{R}^{3\times3}$ as
\begin{align}
A_{\mu\nu}&=\int \text{d}\vecx\,p(\vecx)x_{\mu}x_{\nu}e^{-\lambda^2(\vecw\cdot\vecx)^2}\nonumber\\
&= \int \text{d}\vecx\,p(\vecx)x_{\mu}x_{\nu}e^{-\tfrac{1}{2}(\vecx'-\vecb)\cdot\matB(\vecx'-\vecb)}
\end{align}
where $\vecb = (-w_0/w_1,0)$, $B_{11} = 2\lambda^2w_1^2$, $B_{12}=2\lambda^2w_1w_2$, $B_{22} = 2\lambda^2w_2^2$, and $\vecx'=(x_1,x_2)$.  This integral does not exist if $p(\vecx')$ is uniform (recall that $x_0=1$ is constant) as can be seen by computing the determinant of $\matB$ which is equal to zero for all $w_0$ and $\vecw$. It can be evaluated if we assume $p(\vecx') = \mathcal{N}(\matzero,\Sv)$. We assume $\Sv$ is diagonal with each diagonal element equal to $\sigma_p^2 = 0.25$. In this case $A_{\mu\nu}$ can be rewritten as
\begin{equation}
A_{\mu\nu}= C \int \text{d}\vecx\,x_{\mu}x_{\nu}e^{-\tfrac{1}{2}(\vecx'-\vecb')\cdot\matB'(\vecx'-\vecb')}
\end{equation}
where $\matB' = \matB + \Sv^{-1}$ and $\vecb' = (\matB+\Sv^{-1})^{-1}\matB\vecb$ and $C = e^{-\tfrac{1}{2}\vecb\cdot\matB\vecb +\tfrac{1}{2}\vecb'\cdot\matB'\vecb'}/(2\pi)^{d/2}|\Sv|^{1/2}$. The elements of $\matA$ are
\begin{align*}
A_{00} &=  C(2\pi)^{d/2}|\matB'|^{-1/2}\\
A_{0\mu} &=  C(2\pi)^{d/2}|\matB'|^{-1/2}b_{\mu}'\\
A_{\mu\nu} &= C(2\pi)^{d/2}|\matB'|^{-1/2}(B^{'-1}_{\mu\nu} + b_{\mu}'b_{\nu}').
\end{align*}
We can then write
\begin{equation}
\tilde{\sigma}^2(D_n) = \frac{\lambda^2}{2\pi}\text{tr}(\matA\matF^{-1}).
\end{equation}
This expression is in the same form as the $A$-optimality criterion and analogous to the utility function developed by \cite{schein2007active}. If $d=1$ then $b_1 = -w_0/w$ and $B_{11} = 2\lambda^2w^2$. Extension do larger values of $d$ is straightforward.

\subsection{Computational Note}

Computationally, the most expensive part of the procedure is determining the maximum and minimum utility values. As local optima may exist multiple attempts must be made to locate the global optima. A simulated trial with two covariates, three treatment arms, and a total of 100 recruits can be completed in under ten minutes with an Intel i7 quad core CPU (using any of the utility functions). Uncertainty sampling is considerably faster than the three alternatives and takes a matter of seconds to simulate.

%
%
\newpage
\section{Supplementary Results and Figures}
\label{eq:implementation}
%
%


\begin{figure}[thb!]
\centering
\includegraphics[scale=0.5]{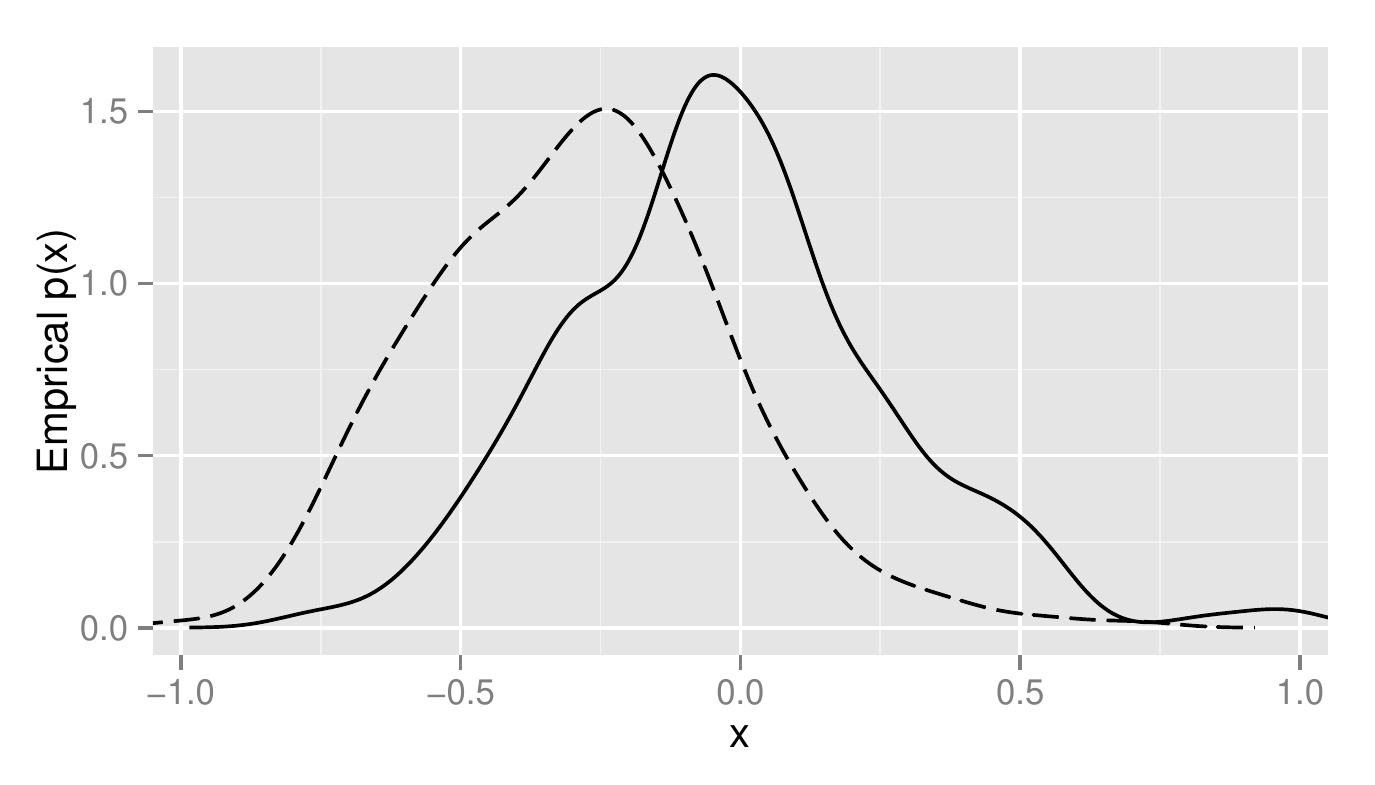}
\caption{Case study: The empirical covariate distribution for the $y=+1$ class (solid line) and the $y=-1$ class (dashed line). A Gaussian kernel with bandwidth $=0.07$ is used to plot the empirical distribution.}
\label{fig:case_study1}
\end{figure}

\begin{figure}[thb]
\centering
\begin{tabular}{c c}
\subfloat{\includegraphics[scale=0.6]{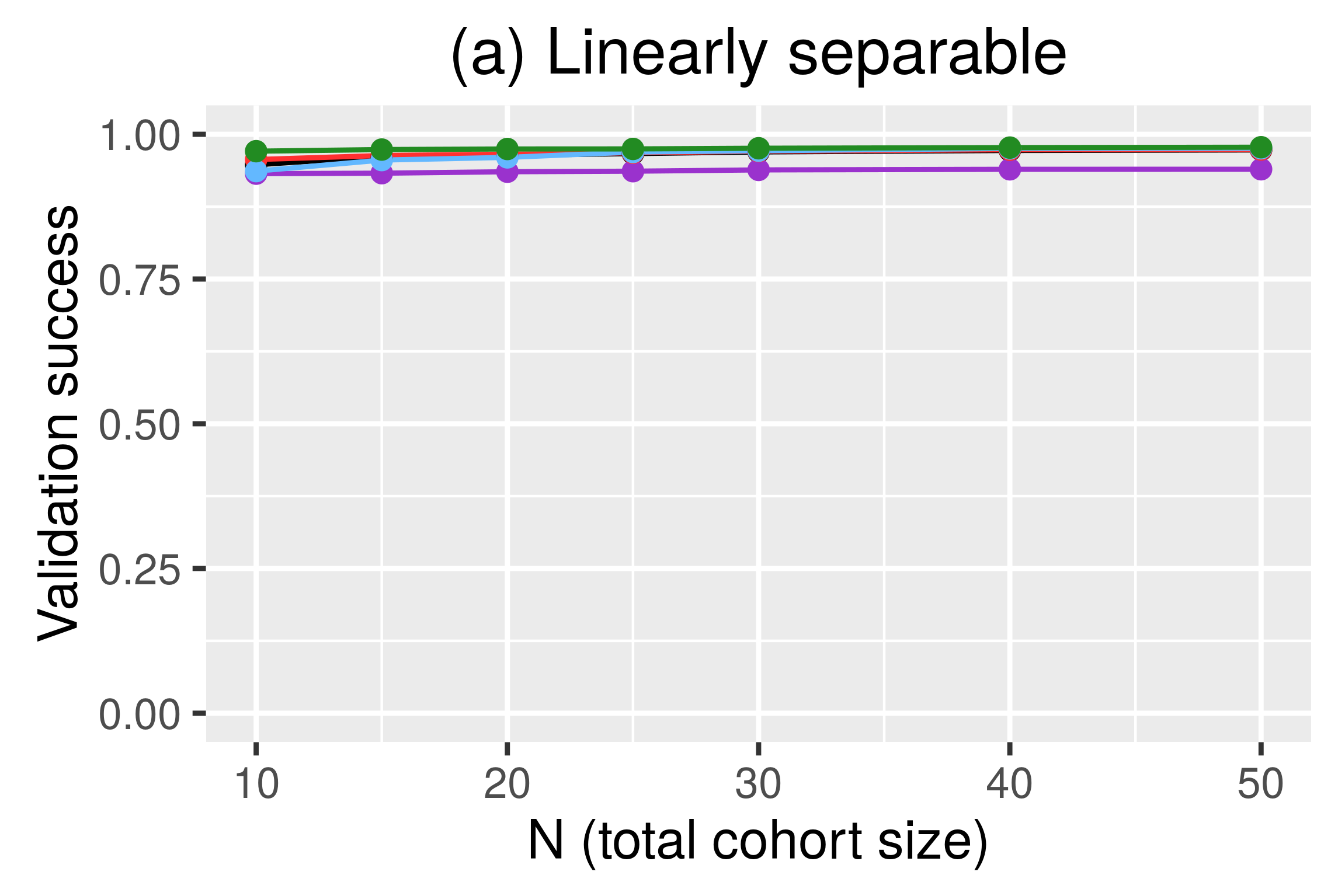}} & \subfloat{\includegraphics[scale = 0.6]{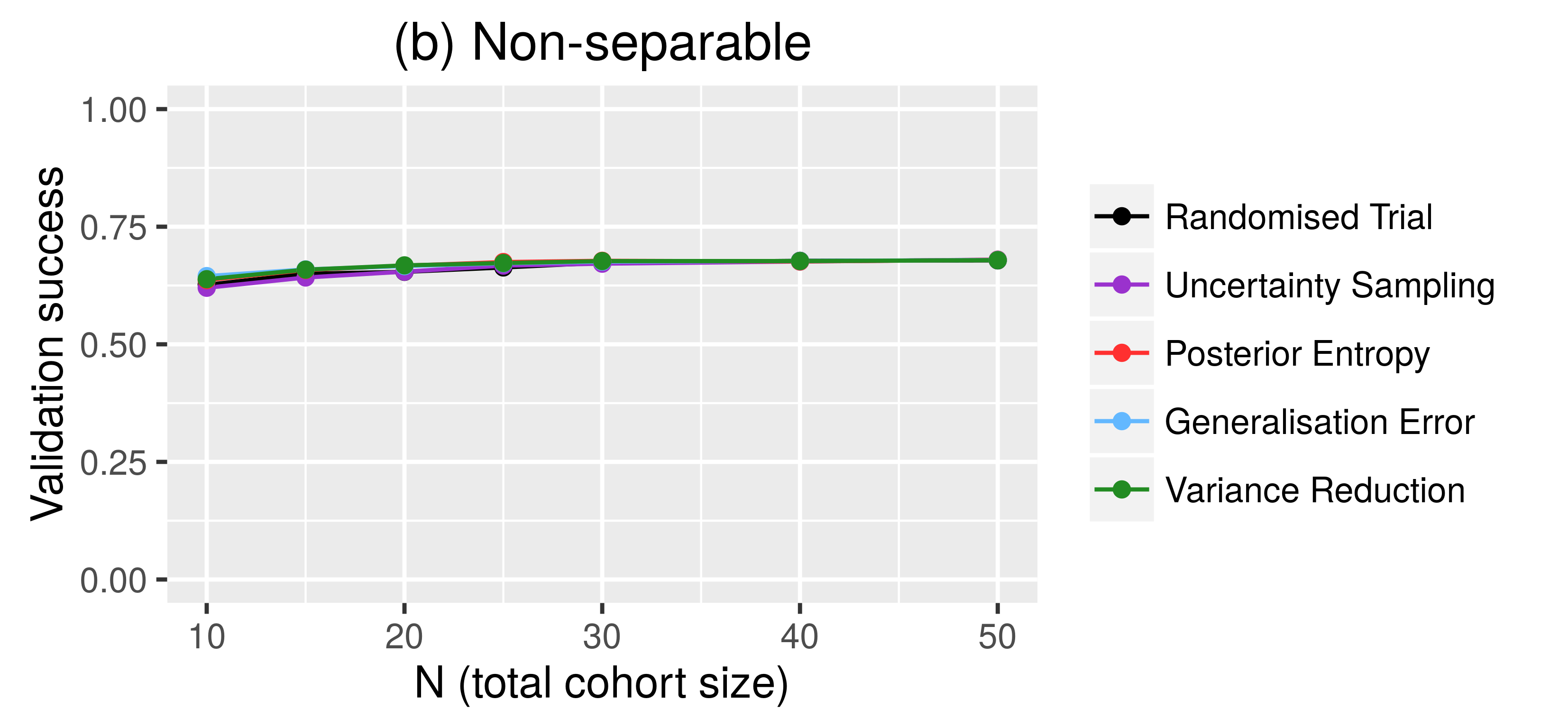}}\\
\end{tabular}
\caption{Simulation study 1: success rate in validation cohorts for both the linearly separable and non-seperable datasets. All data were obtained by averaging over 500 simulations.}

\label{fig:sim2:1}
\end{figure}

\begin{figure}[tb!]
\centering
\begin{tabular}{c c}
\subfloat[Predictive density $p(y|x,w,w_0)$]{\includegraphics[scale=0.5]{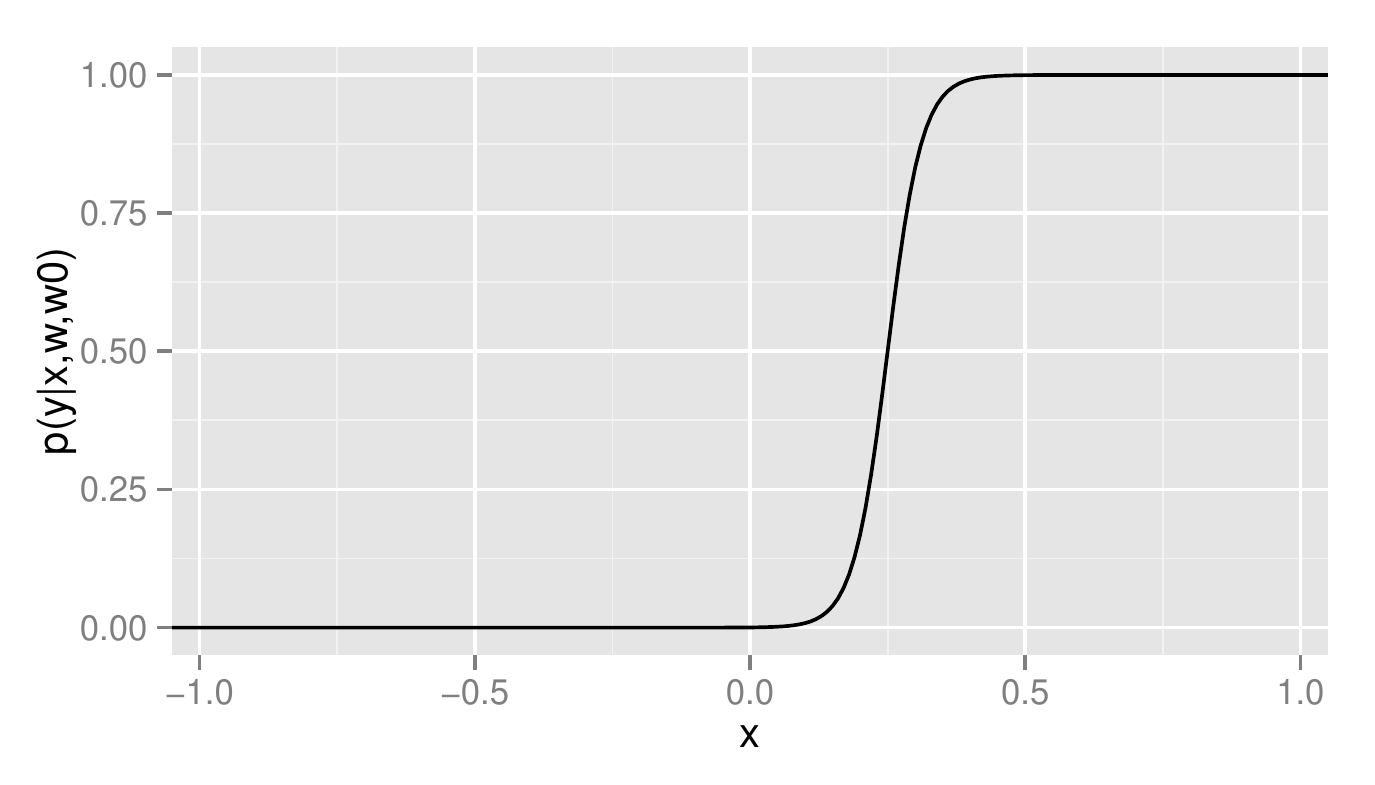}} & \subfloat[RCT]{\includegraphics[scale = 0.5]{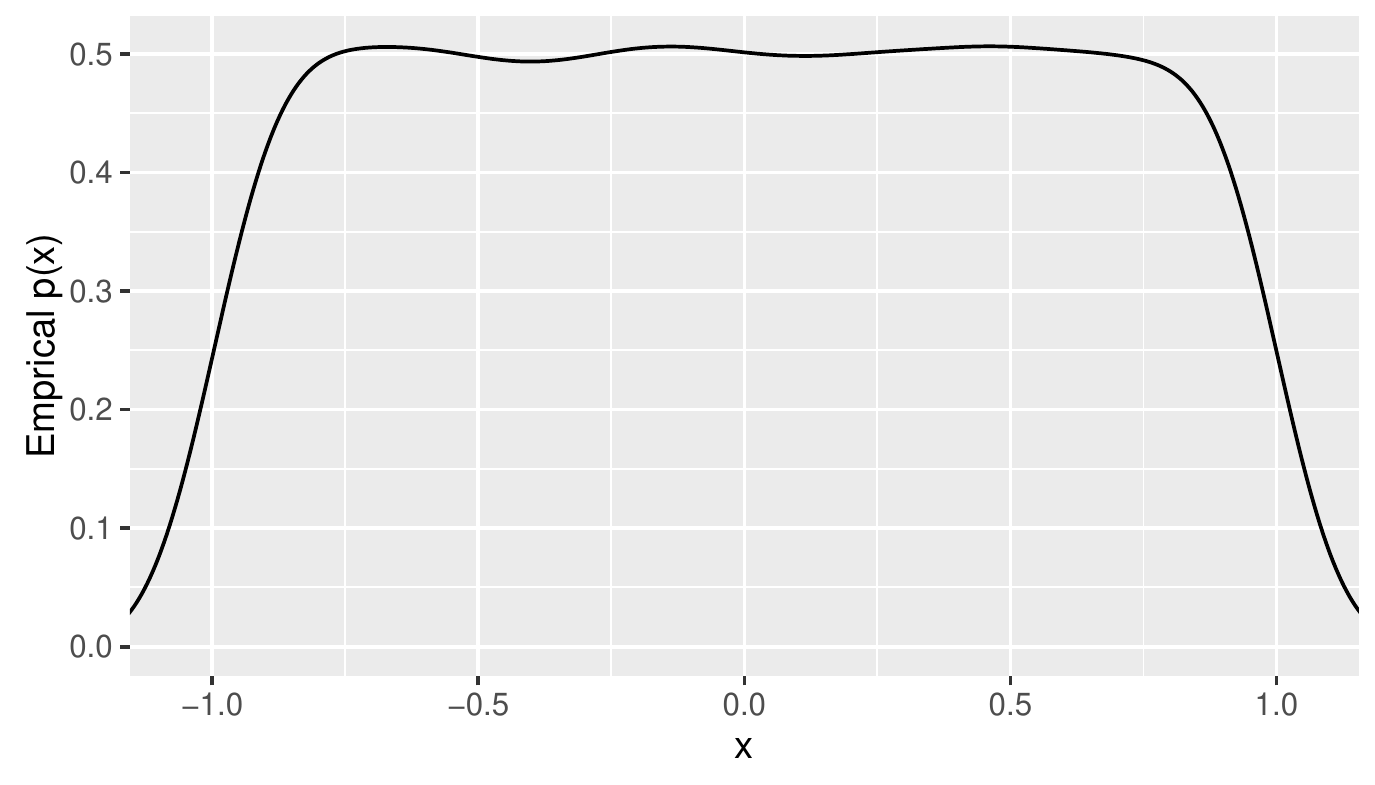}}\\
\subfloat[Uncertainty Sampling]{\includegraphics[scale=0.5]{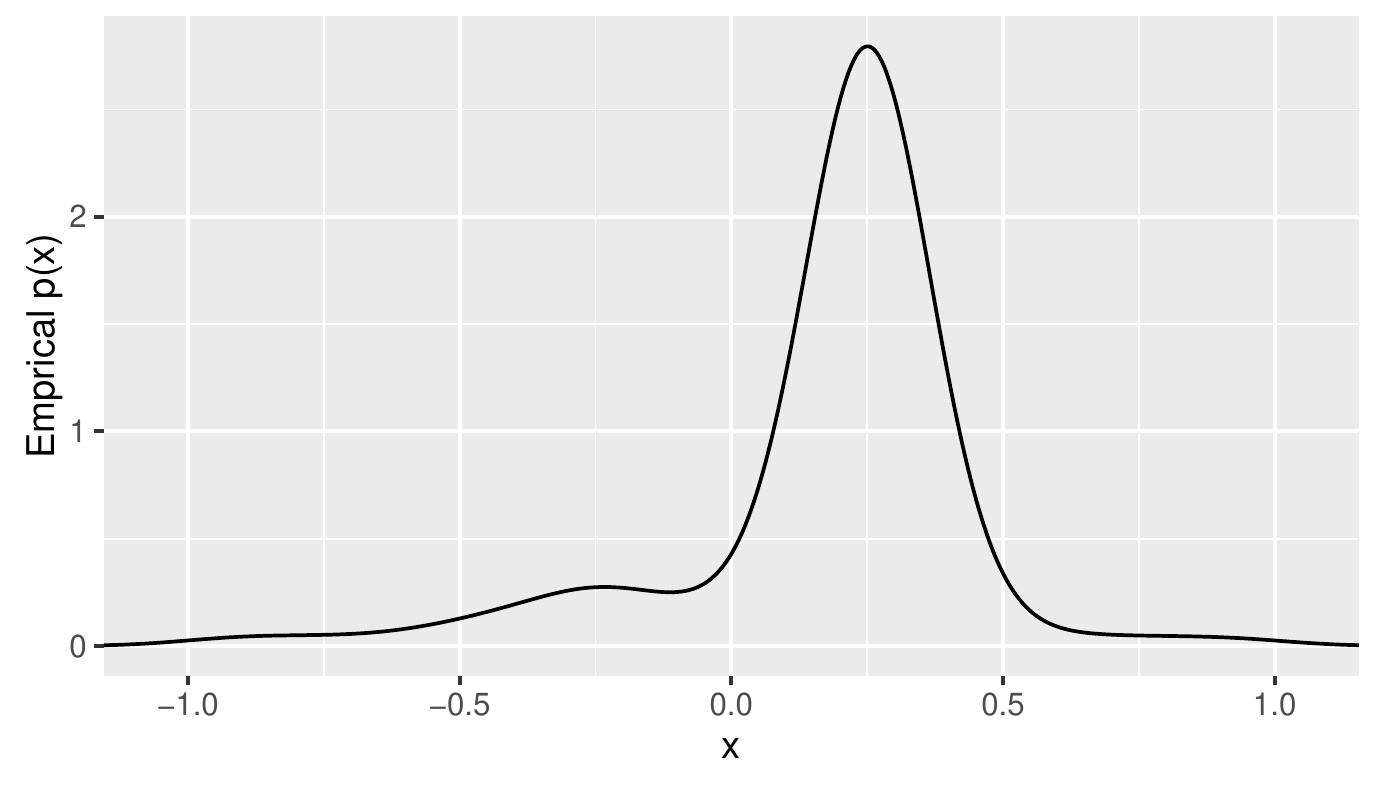}} & \subfloat[Posterior Entropy]{\includegraphics[scale = 0.5]{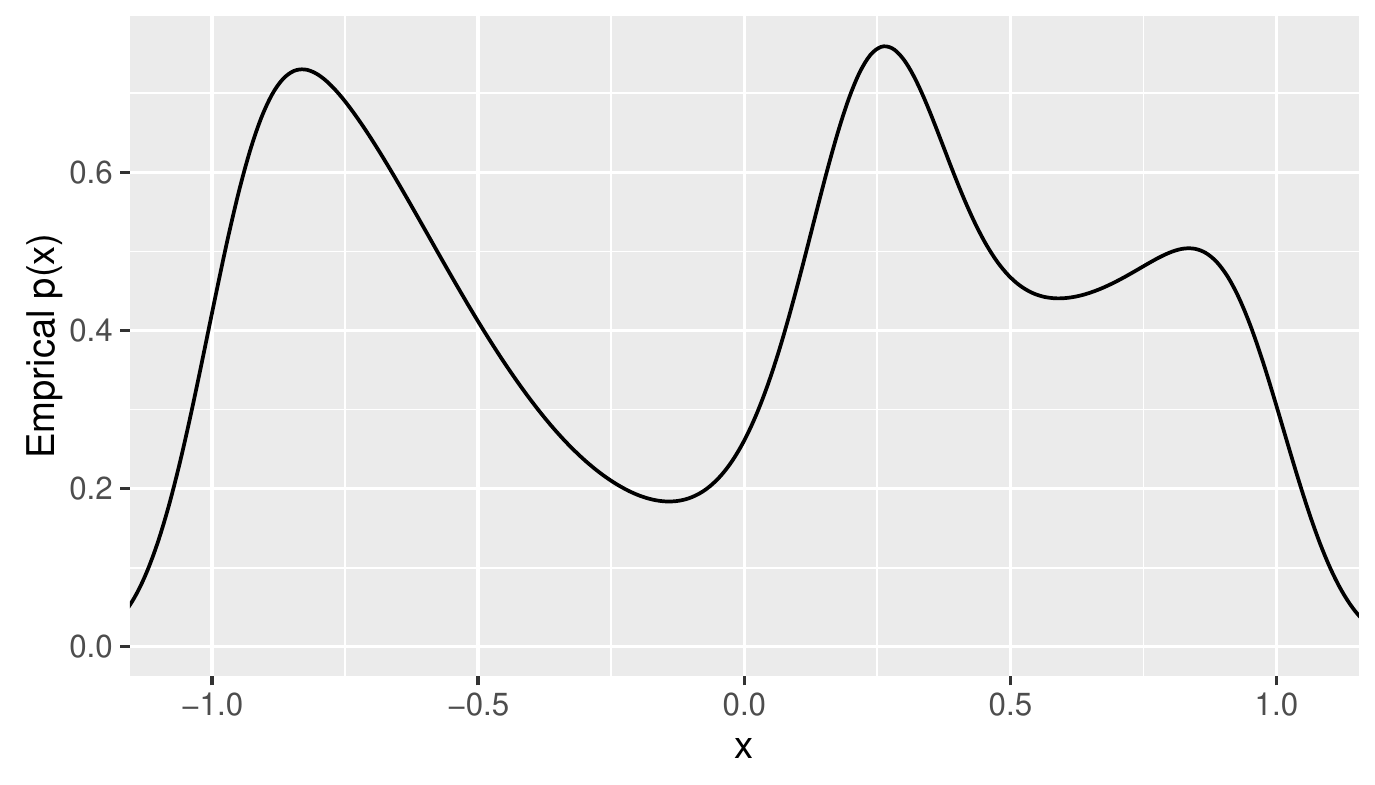}}\\
\subfloat[Variance Reduction]{\includegraphics[scale=0.5]{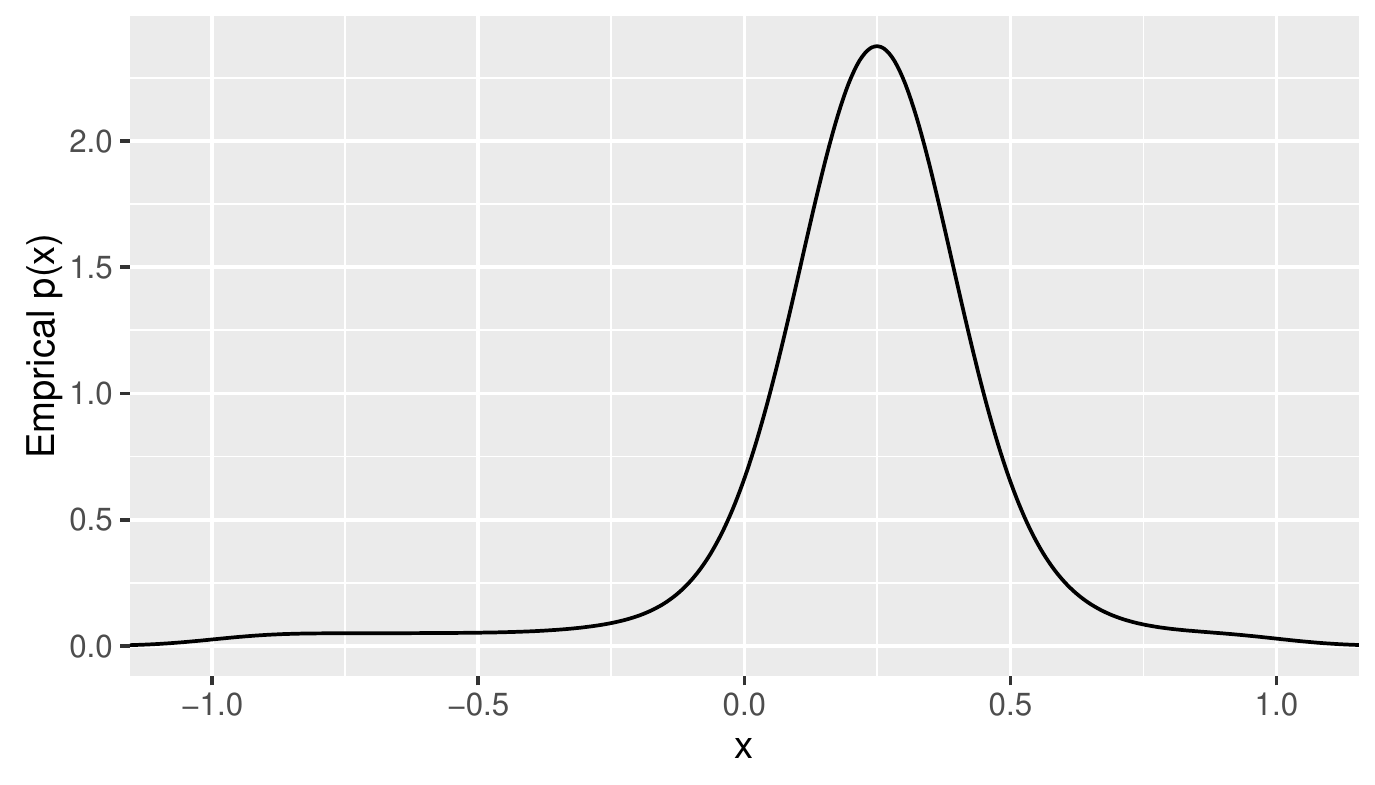}} & \subfloat[Generalisation Error]{\includegraphics[scale = 0.5]{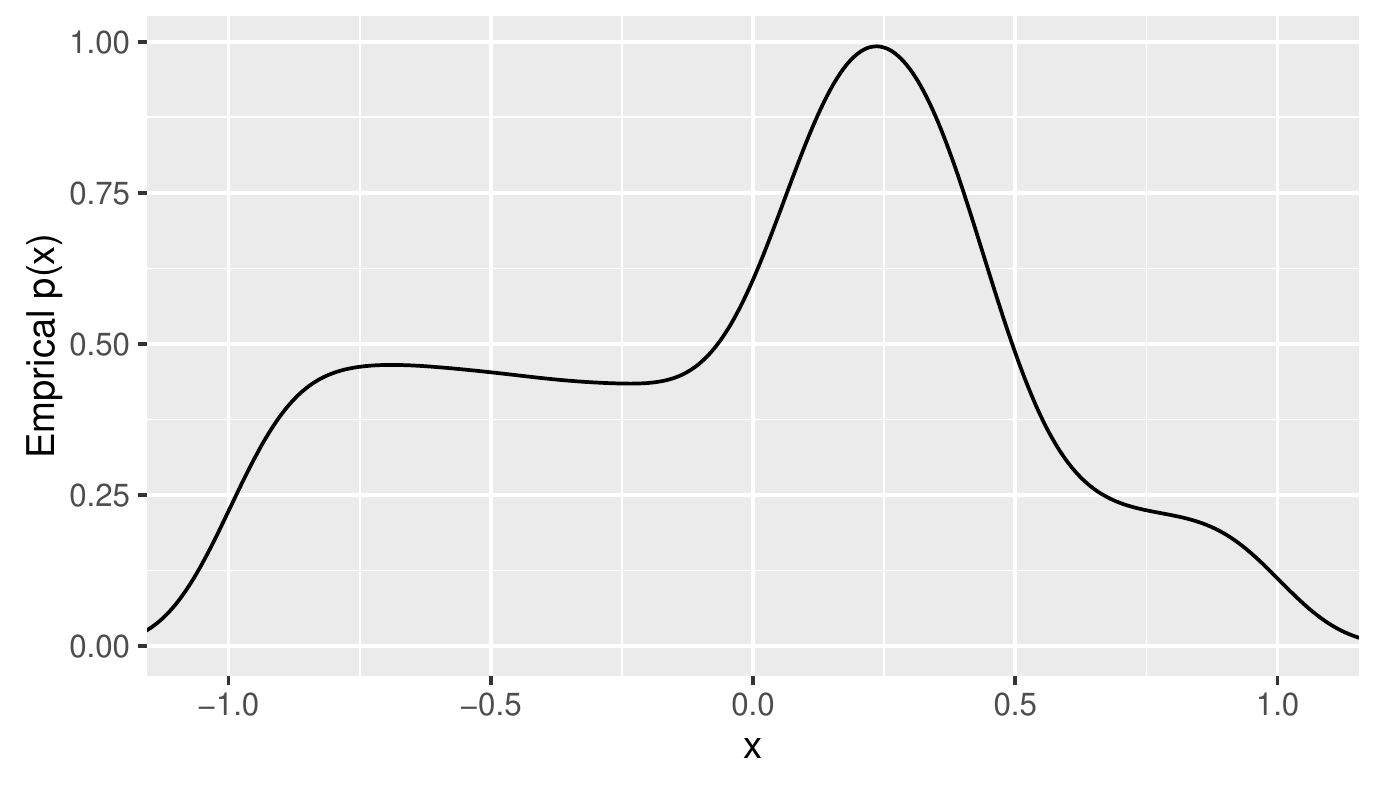}}\\
\end{tabular}
\caption{Simulation study 1: In (a) is $p(y=+1|x,w,w_0)$ as a function of $x$ where $w=32$ and $w_0=-8$. The decision boundary is located at $x=0.25$. In (b-f) are empirical cohort distributions obtained by averaging over 500 simulations in the linearly separable case. A Gaussian kernel with bandwidth $=0.1$ was used.}
\label{fig:sup:sim1:2}
\end{figure}

\begin{figure}[tb!]
\centering
\begin{tabular}{c c}
\subfloat[Theoretical distribution]{\includegraphics[scale=0.5]{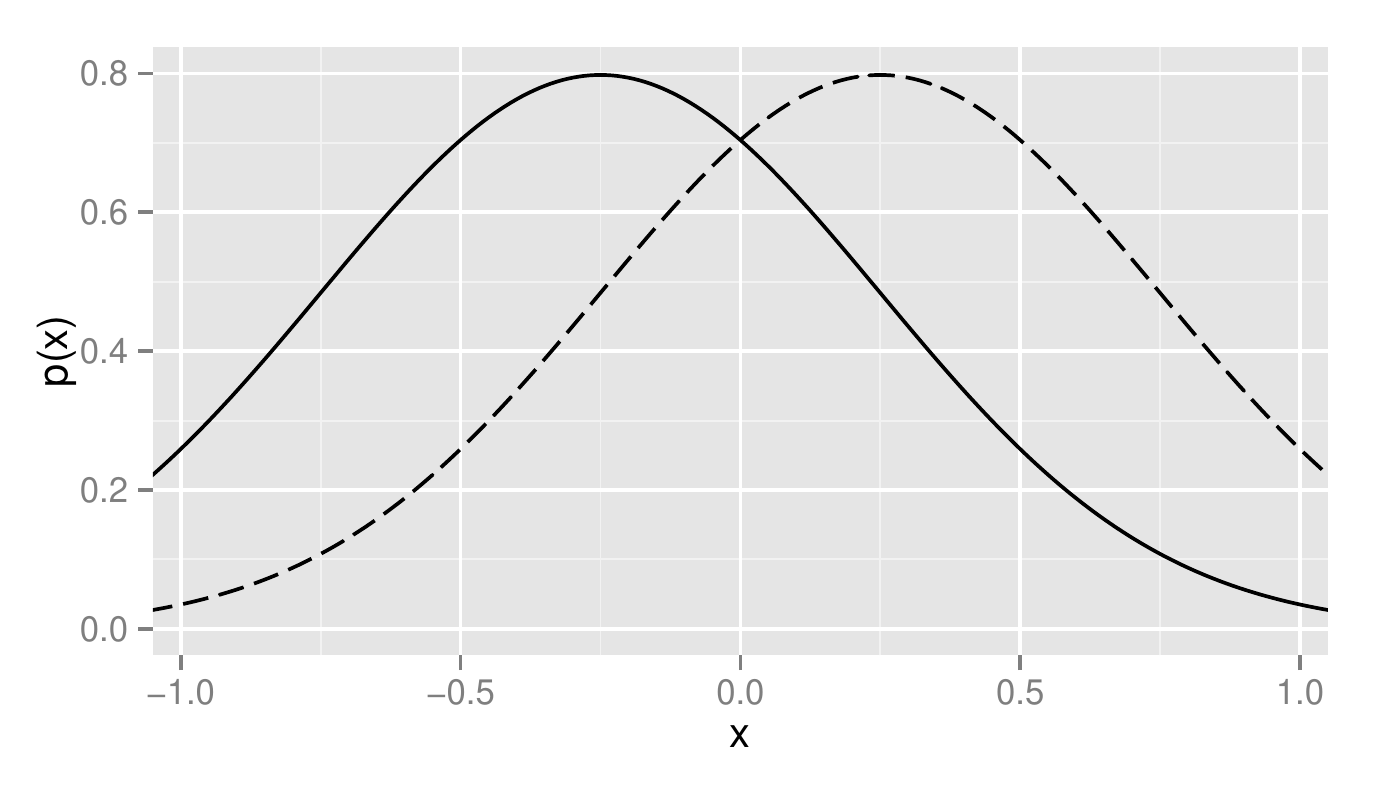}} & \subfloat[RCT]{\includegraphics[scale = 0.5]{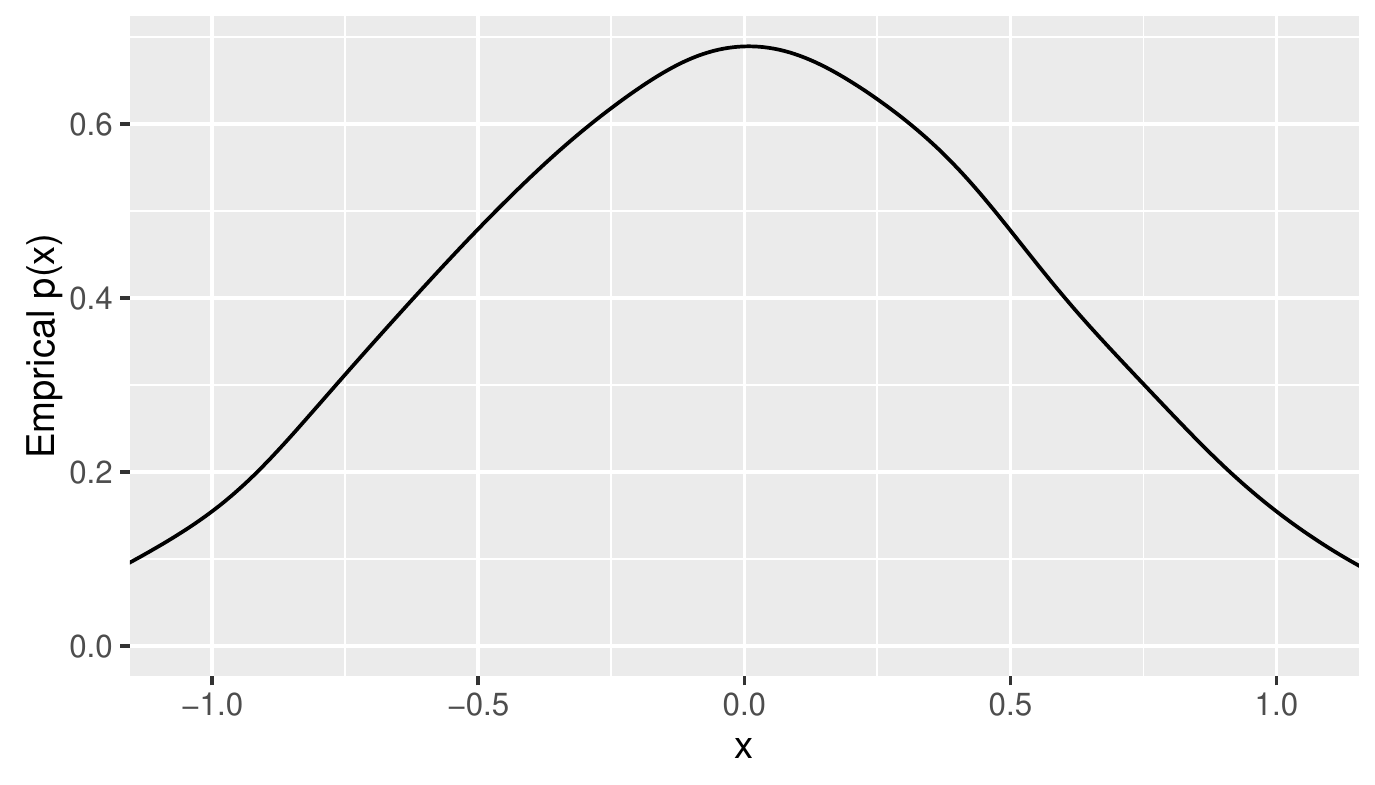}}\\
\subfloat[Uncertainty Sampling]{\includegraphics[scale=0.5]{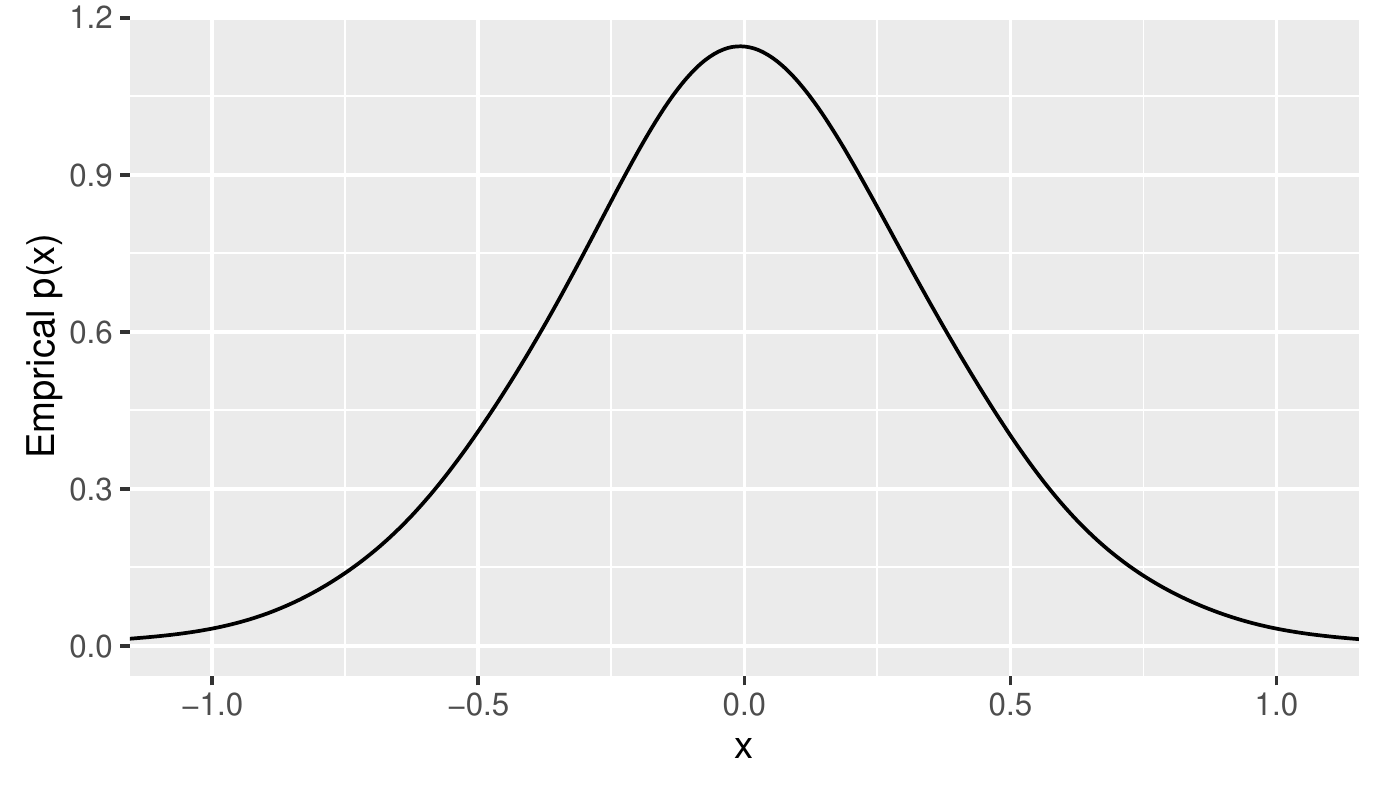}} & \subfloat[Posterior Entropy]{\includegraphics[scale = 0.5]{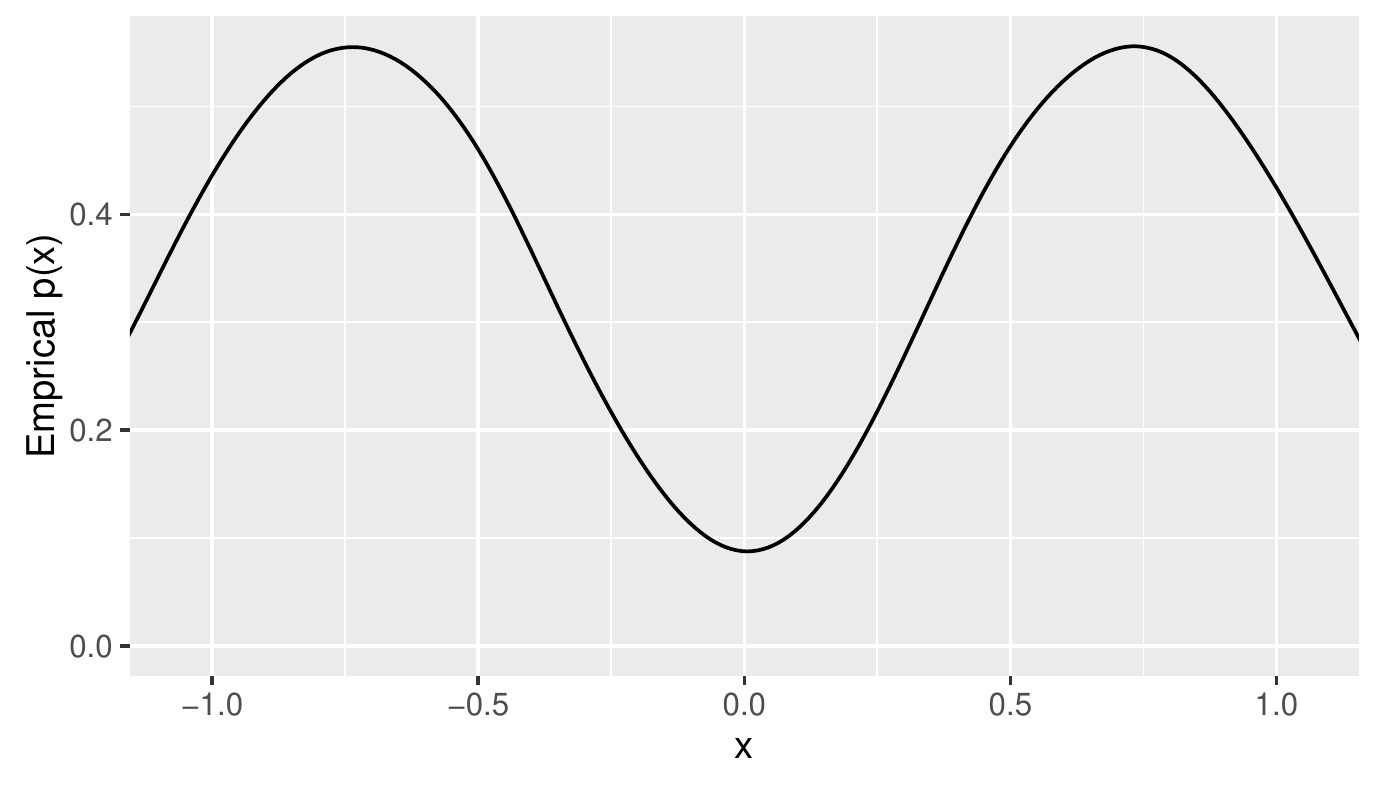}}\\
\subfloat[Variance Reduction]{\includegraphics[scale=0.5]{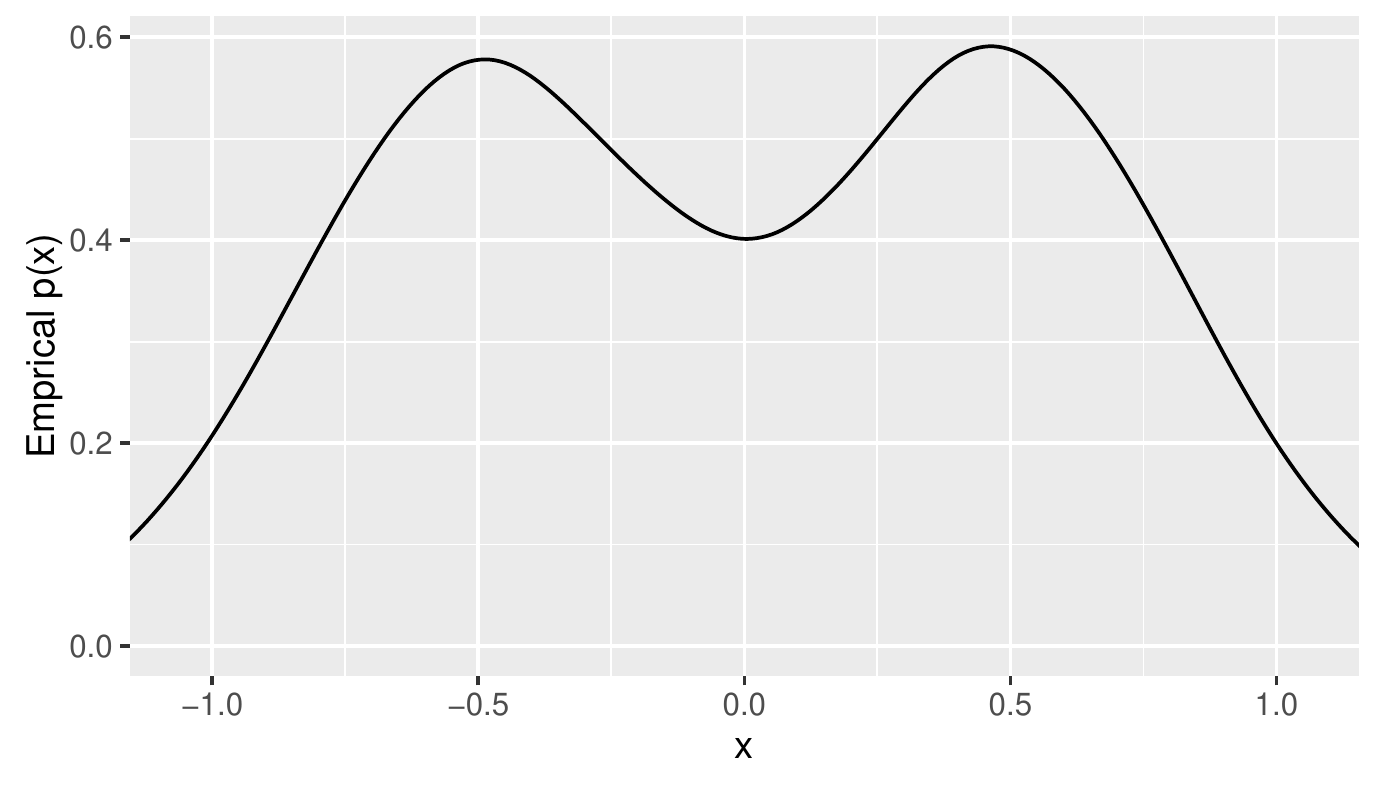}} & \subfloat[Generalisation Error]{\includegraphics[scale = 0.5]{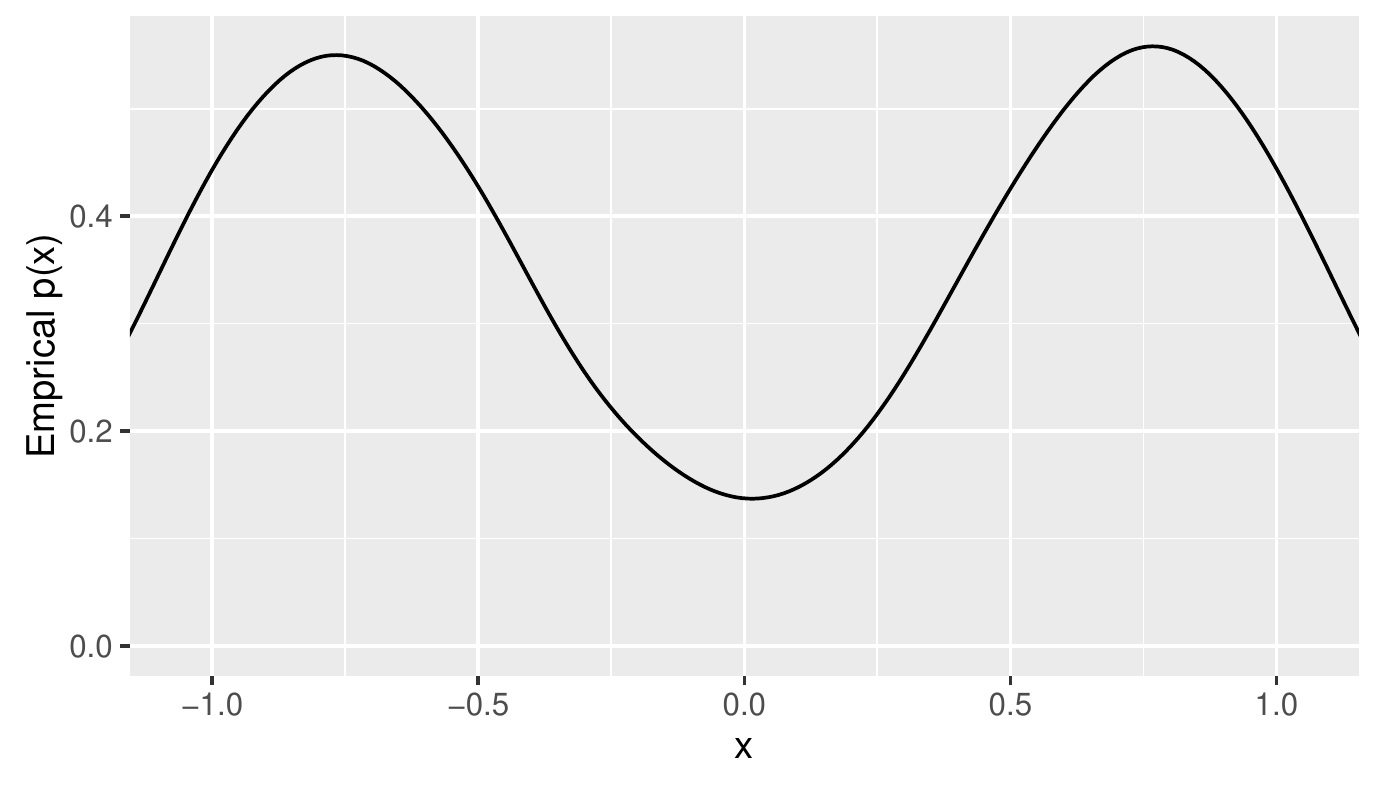}}\\
\end{tabular}
\caption{Simulation study 1: In (a) is distribution of covariates that was used to generate the simulated data in the non-separable case. Individuals with $y=+1$ were generated from a Gaussian centred on $x=-0.25$ (solid line) and individuals with $y=-1$ were generated from a Gaussian at $x=+0.25$ (dashed line). In (b-f) are empirical cohort distributions obtained by averaging over 500 simulations in the non-separable case. A Gaussian kernel with bandwidth $=0.1$ was used.}
\label{fig:sup:sim1:3}
\end{figure}


\begin{figure}[tb]
\centering
\begin{tabular}{c c}
\subfloat{\includegraphics[scale=0.6]{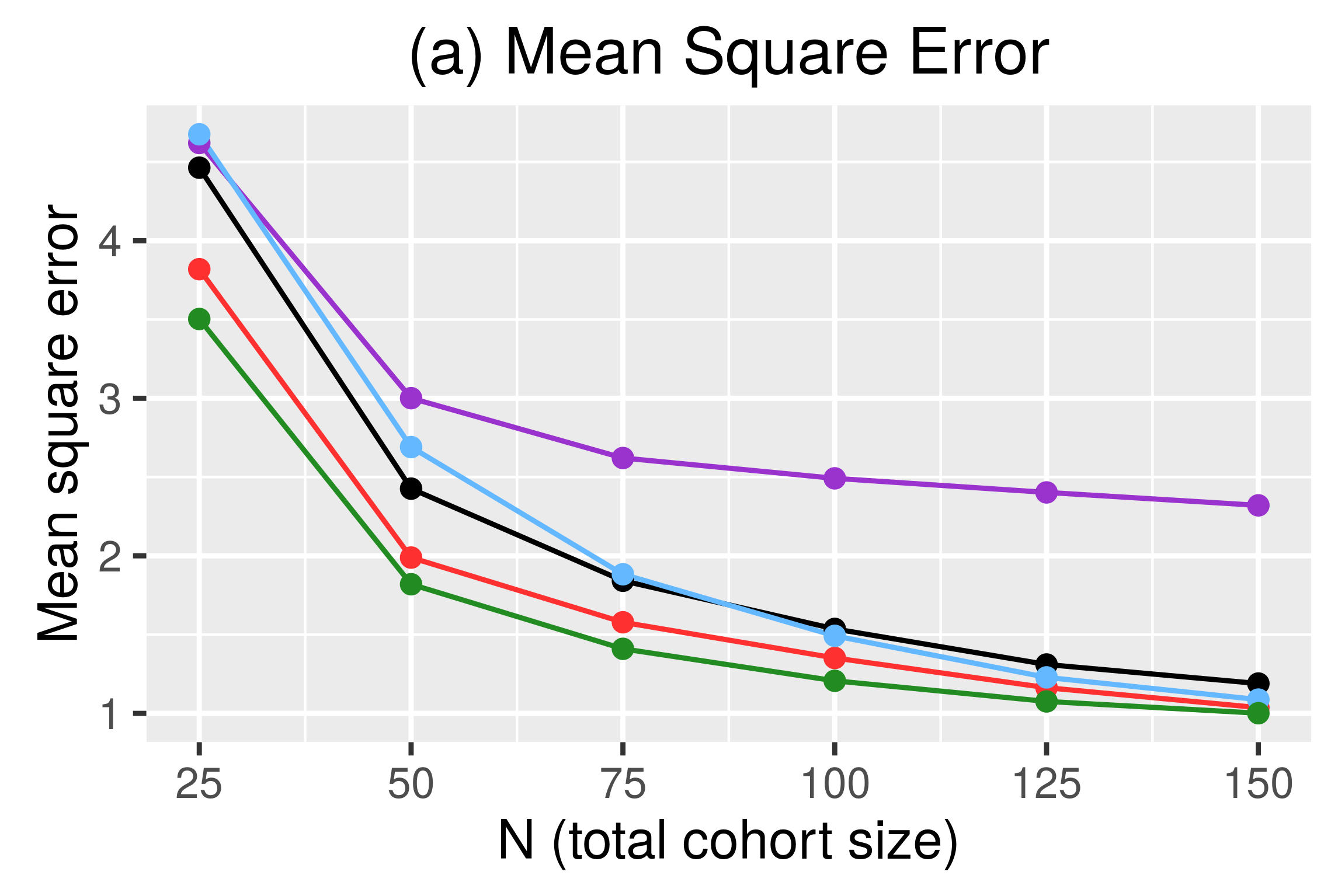}} & \subfloat{\includegraphics[scale = 0.6]{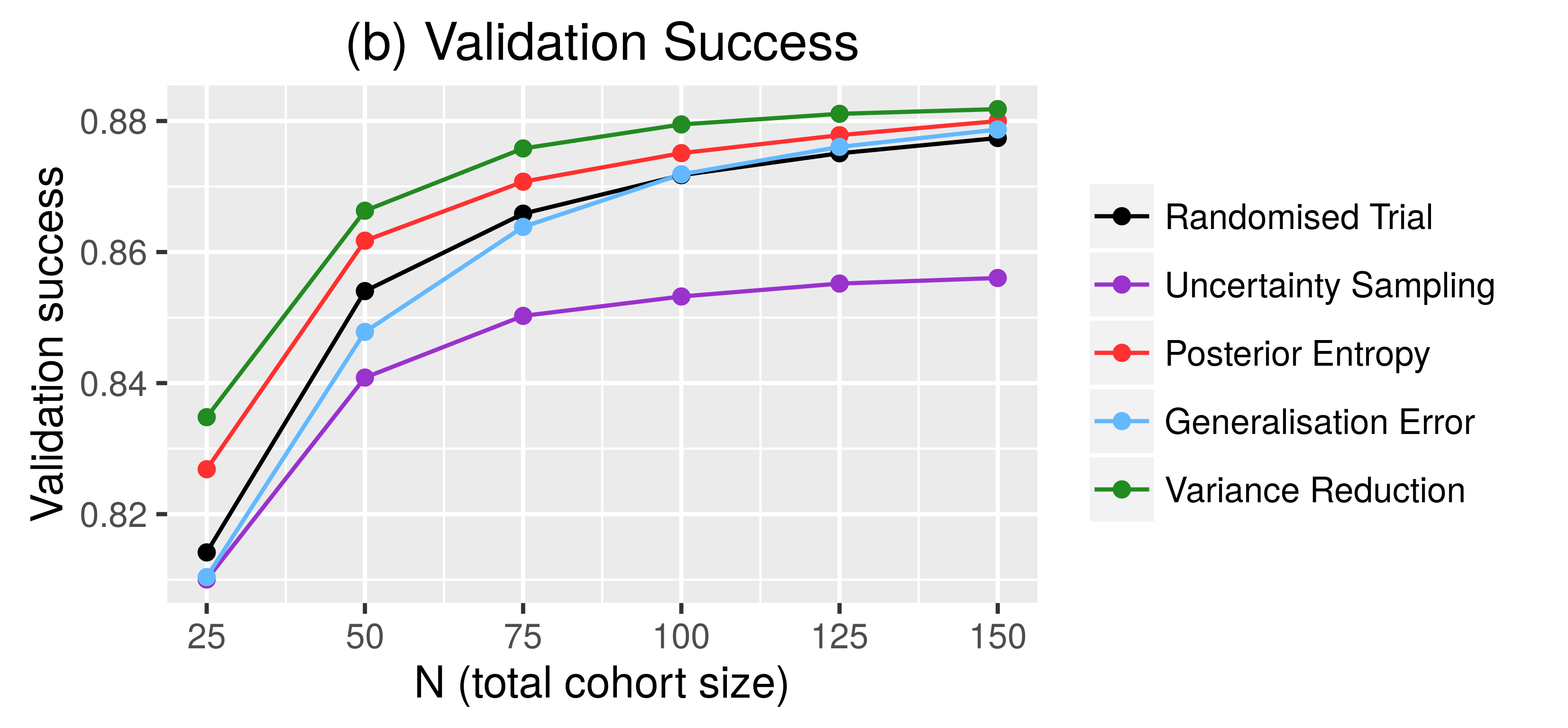}}\\
\end{tabular}
\caption{Simulation study 2: three treatment arms, selective recruitment, and adaptive treatment allocation. In (a) is the mean square error between the true and inferred model parameters (averaged over all components of $\vecw$ and $w_0$). Panel (b) is the proportion of correct predictions on a validation cohort. All data have were obtained by averaging over 500 simulations.}

\label{fig:sim2:1}
\end{figure}

\begin{figure}[tb!]
\centering
\begin{tabular}{c c}
\subfloat{\includegraphics[scale=0.6]{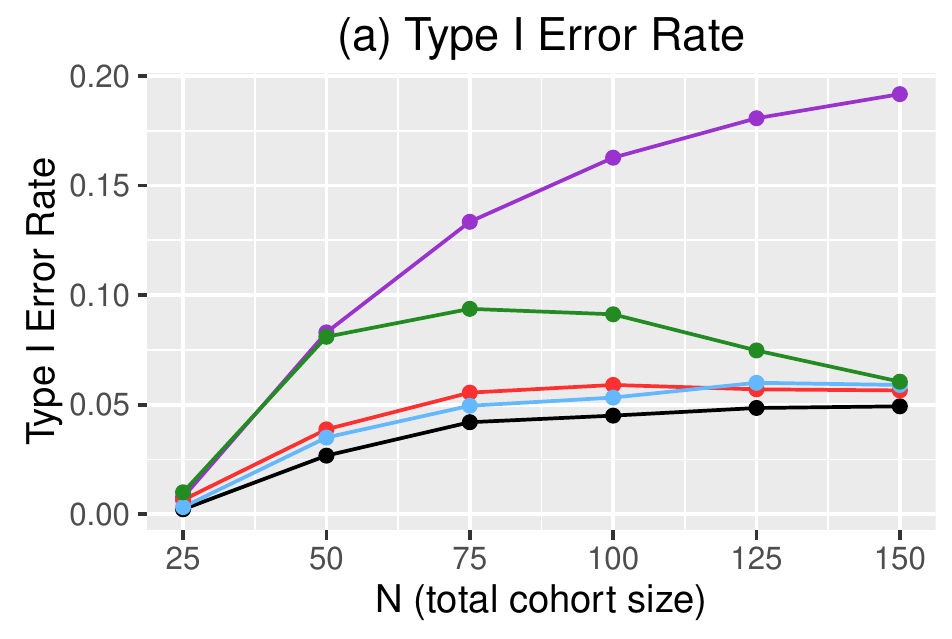}} & \subfloat{\includegraphics[scale = 0.6]{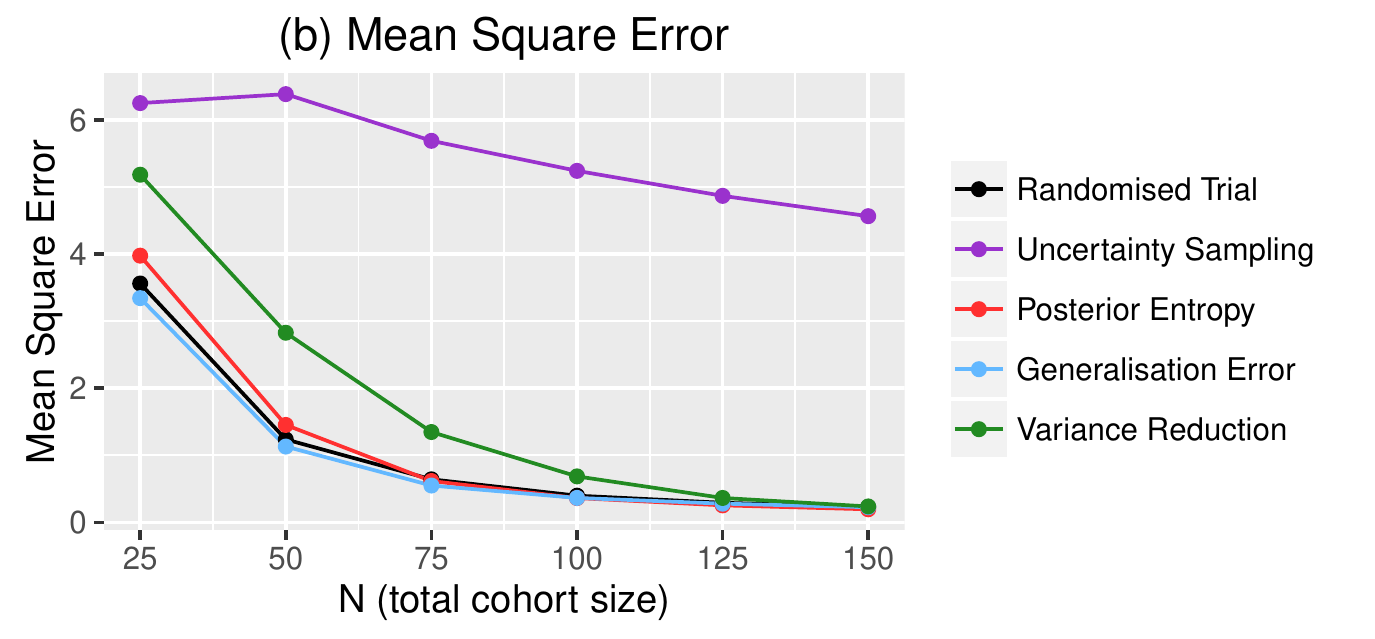}}\\
\end{tabular}
\caption{Simulation study 2: results from further simulations using the null model where $w_0=0$ and $\vecw=0$ for all treatment arms. The mean square error is between the inferred model parameters and the true values (which equal zero in this case). The estimates are asymptotically unbiased and converge to zero. The Type I error rate is roughly 5\% as expected. The Uncertainty Sampling method is an exception and appears to converge very slowly to the true parameter values and has an inflated error rate. The Variance Reduction method is somewhat slower to converge and also has an elevated Type I error rate but this stabilises for larger cohort sizes.}
\label{fig:sup:sim6}
\end{figure}

\begin{figure}[tb]
\centering
\begin{tabular}{c c}
\subfloat{\includegraphics[scale=0.6]{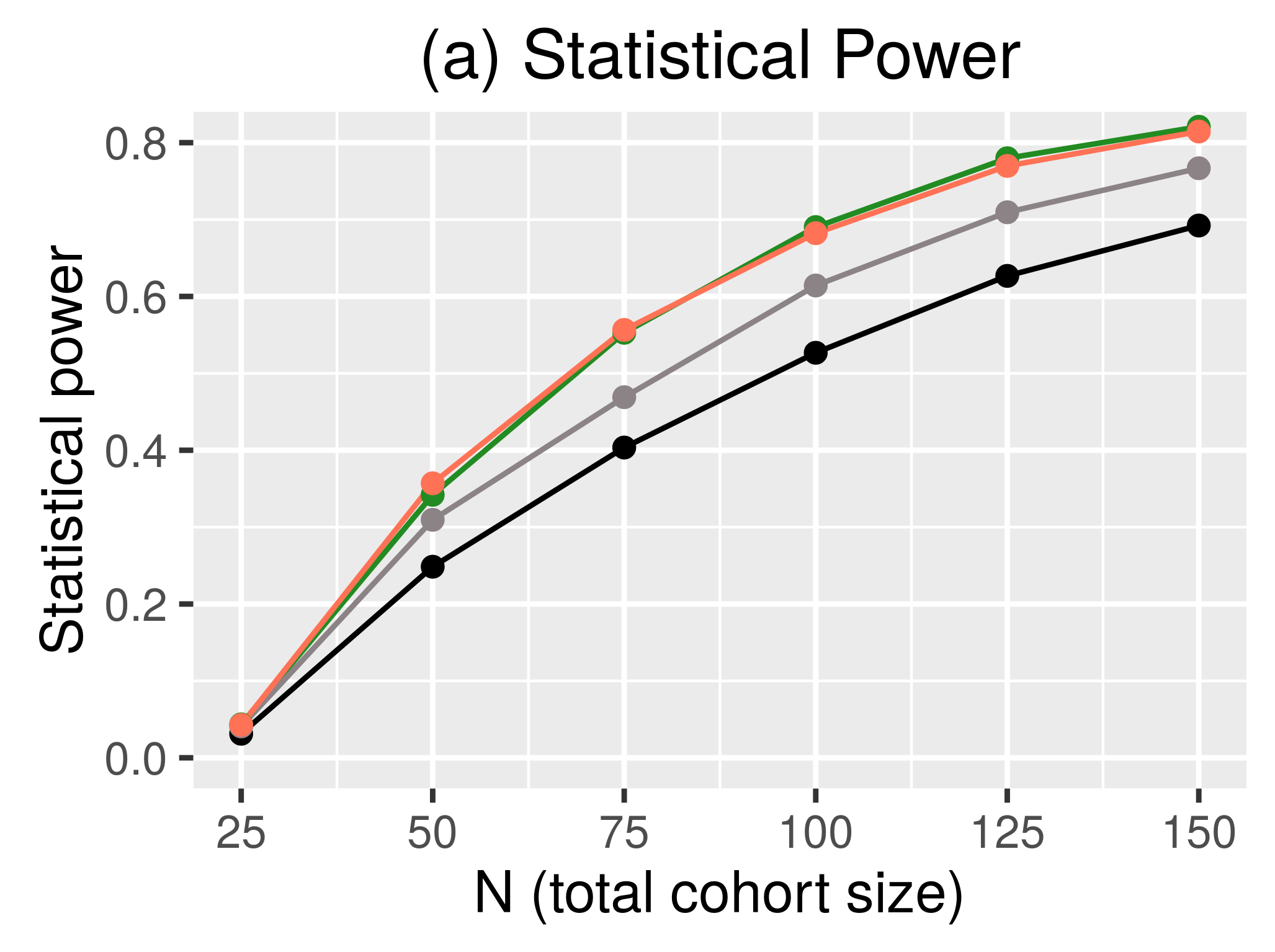}} & \subfloat{\includegraphics[scale = 0.6]{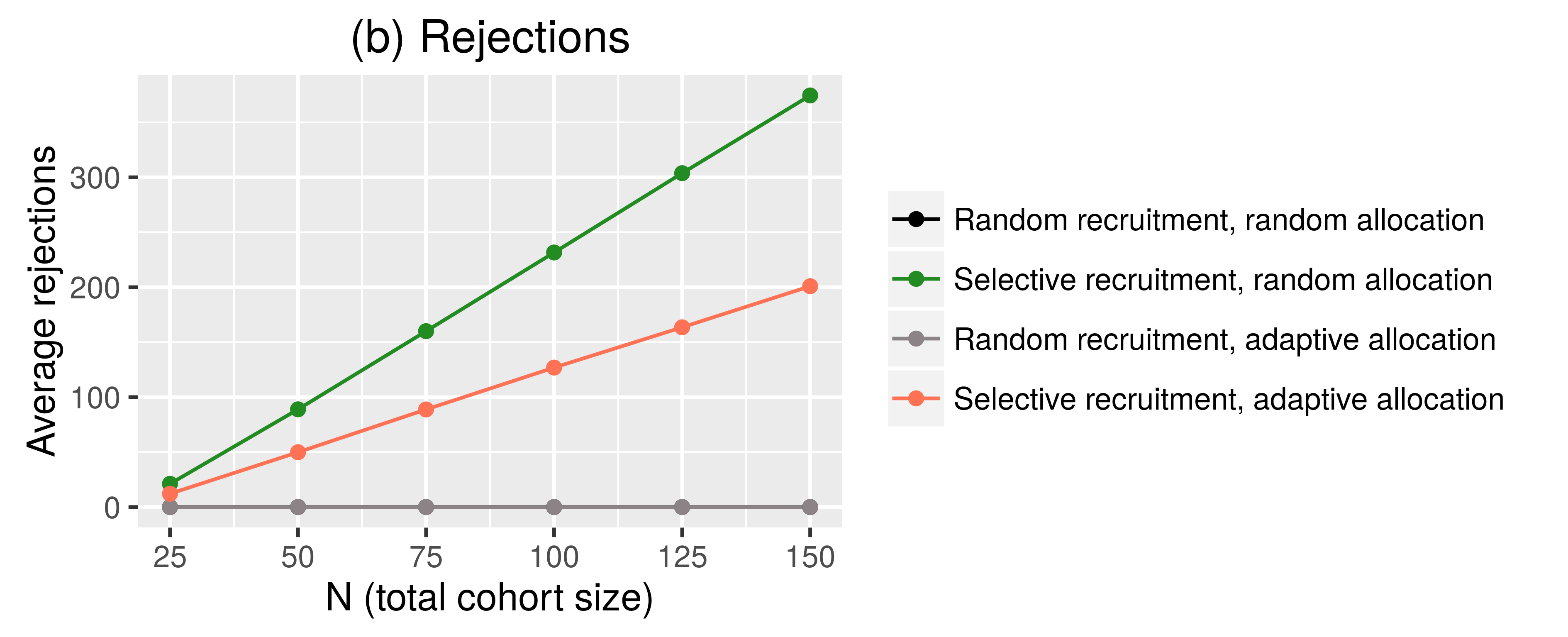}}\\
\end{tabular}
\caption{Simulation study 3: results from further simulations using the same configuration as Simulation Studies 2 and 3 in the main text. In this plot the statistical power and number of rejections are plotted using the variance reduction method only. Three types of adaptive designs (using different combinations of selective recruitment and adaptive allocation) are compared to a randomised design. Random recruitment means that all individuals are recruited onto the trial. Selective recruitment achieves similar power with or without adaptive allocation although more individuals are rejected with random allocation. This is due to the protocol that a treatment arm is selected before recruitment (see Section 2.3 of the main text). Uninformative individuals are always filtered during recruitment regardless of how the treatment arm was chosen. Adaptive allocation will tend to select treatment arms that result in higher recruitment probabilities hence fewer individuals are rejected in that case.}

\label{fig:sim2:1}
\end{figure}

\end{document}